\input ./style/arxiv-general.cfg
\documentclass[aos,MSNbibl,seceqn,dvips]{arximspdf}
\makeatletter
   \@ifpackageloaded{graphicx}{}{\usepackage{graphicx}}
\makeatother
\usepackage[algoruled,linesnumbered]{algorithm2e}

\doi{10.1214/15-AOS1310}
\volume{43}
\issue{3}
\pubyear{2015}
\firstpage{1300}
\lastpage{1322}
\docsubty{FLA}

\makeatletter
\let\my@algocf@latexcaption\algocf@latexcaption
\let\my@addcontentsline\addcontentsline
\long\def\algocf@latexcaption#1[#2]#3{%
\def\addcontentsline##1##2##3{}%
\my@algocf@latexcaption{#1}[#2]{#3}%
\global\let\addcontentsline\my@addcontentsline%
}
\newcommand{\eqref}[1]{(\ref{#1})}
\newtheorem{theorem}{Theorem}[section]
\newproclaim{fact}[theorem]{Fact}
\newtheorem{corollary}[theorem]{Corollary}

\newtheorem{lemma}[theorem]{Lemma}
\newproclaim{definition}[theorem]{Definition}
\newtheorem{proposition}[theorem]{Proposition}
\newproclaim{assumption}[theorem]{Assumption}

\newproclaim{rem}[theorem]{Remark}

\newtheorem{conjecture}[theorem]{Conjecture}

\makeatother

\begin{document}
\begin{frontmatter}

\title{Do semidefinite relaxations solve sparse PCA up to the
information limit?}
\runtitle{Do semidefinite relaxations solve sparse PCA?}

\begin{aug}
\author[A]{\fnms{Robert} \snm{Krauthgamer}\ead[label=e1]{robert.krauthgamer@weizmann.ac.il}\thanksref{T1,m1}},
\author[A]{\fnms{Boaz} \snm{Nadler}\ead[label=e2]{boaz.nadler@weizmann.ac.il}\thanksref{T2,m1}}
\and
\author[B]{\fnms{Dan} \snm{Vilenchik}\corref{}\ead[label=e3]{vilenchi@bgu.ac.il}\thanksref{T3,m2}}

\thankstext{T1}{Supported in part by a
US--Israel BSF Grant \#2010418, ISF Grant \#897/13 and by the Citi Foundation.}
\address[A]{R.~Krauthgamer\\
B.~Nadler\\
Weizmann Institute of Science\\
Rehovot\\
Israel\\
\printead{e1}\\
\phantom{E-mail:\ }\printead*{e2}}
\thankstext{T2}{Supported in part by the Citi Foundation.}

\address[B]{D.~Vilenchik\\
Ben-Gurion University\\
Beer-Sheva\\
Israel\\
\printead{e3}}
\thankstext{T3}{Research done in part while a postdoctoral fellow at
the Weizmann Institute.}

\runauthor{R.~Krauthgamer, B.~Nadler and D.~Vilenchik}

\affiliation{Weizmann Institute of Science\thanksmark{m1} and
Ben-Gurion University\thanksmark{m2}}
\end{aug}

%
\received{\smonth{9} \syear{2014}}
%
\revised{\smonth{1} \syear{2015}}

%
\begin{abstract}
Estimating the leading principal components
of data, assuming they are sparse,
is a central task in modern high-dimensional statistics.
Many algorithms were developed for this sparse PCA problem, from
simple diagonal thresholding to sophisticated
semidefinite programming (SDP) methods.
A~key theoretical question is under what conditions
can such \mbox{algorithms} recover the sparse principal components?
We study this question for a single-spike model with an $\ell_0$-sparse
eigenvector, in the asymptotic regime as dimension $p$ and sample size
$n$ both tend to infinity.
Amini and Wainwright [\textit{Ann. Statist.}
\textbf{37} (2009) 2877--2921] proved that for sparsity levels
$k\geq\Omega(n/\log p)$,
no algorithm, efficient or not, can reliably recover the sparse eigenvector.
In contrast, for $k\leq O(\sqrt{n/\log p})$,
diagonal thresholding is consistent.
It was further conjectured that an SDP approach may close this gap
between computational and information limits.
We prove that when $k \geq\Omega(\sqrt{n})$, the proposed SDP
approach, at least in its standard usage,
cannot recover the sparse spike.
In fact, we conjecture that in the single-spike model,
no computationally-efficient algorithm can recover a spike
of $\ell_0$-sparsity $k\geq\Omega(\sqrt{n})$.
Finally, we present empirical results suggesting that up to sparsity levels
$k=O(\sqrt{n})$, recovery is possible by
a simple covariance thresholding algorithm.
\end{abstract}

%
\end{frontmatter}

\section{Introduction}\label{sec1}
Principal components analysis (PCA) is a popular technique for
dimension reduction
that has a wide range of applications involving multivariate data, in
both science and engineering; see, for example, \cite{Anderson84,JoliffePCA}.
The first principal component (PC) of a $p$-dimensional random variable
$\mathbf{x}=(x_1,\ldots,x_p)$
is the direction in which the variance of $\mathbf{x}$ is maximal,
or equivalently, the leading eigenvector of its population covariance matrix
$\Sigma=\operatorname{\mathbb E}[(\mathbf{x}-\mu)(\mathbf{x}-\mu
)^T]$ where $\mu=\operatorname{\mathbb E}
[\mathbf{x}]$. In
practice, one
typically does not have explicit access to $\Sigma$, but rather is given
$n$ samples from $\mathbf{x}$, from which one computes the sample covariance
matrix ${\hat{\Sigma}}$ and its leading eigenvectors.

In contemporary applications where variables are plentiful (large $p$)
but samples are relatively scarce (small $n$), PCA
suffers from two major limitations:
(1) the principal components are typically a linear combination of
all variables, which hinders their interpretation and subsequent use,
and (2) while PCA is consistent in the classical setting ($p$ is fixed
and $n\to\infty$) \cite{Muirhead1982,Anderson84}, it is generally
inconsistent in high-dimensions. Indeed, as shown, for example, in
 \cite
{Johnstone01,DPaul07,Nadler08,BickelLevina06,Johnstone.Lu2009Consistency},
when $p$ is comparable to, or significantly larger than $n$, the sample
covariance matrix ${\hat{\Sigma}}$ may be a poor approximation to the
population's covariance matrix $\Sigma$, and its leading eigenvectors
may be far from the population's principal components.


%
%

To address the first drawback, one can consider a \emph{sparse} PCA problem,
in which for some appropriate parameter $k$, we search for a direction
with at most $k$ nonzero coefficients and with maximal variance.
Formally, the $\ell_0$-sparse PCA problem is defined by
%
\begin{equation}
\label{eq:LP} {\mathcal{L}}_0(\Sigma) = \operatorname{argmax} \bigl
\{ \mathbf{w}^T \Sigma\mathbf{w}\dvtx \|\mathbf{w}\|_2 =
1, \| \mathbf{w}\|_0 \leq k \bigr\}.
\end{equation}
We note that other notions of sparsity were considered in the
literature, for example, a population covariance matrix that has only a
few large eigenvalues, whose corresponding eigenvectors are sparse in
$\ell_q$-norm for $q\in(0,2)$ \cite
{Johnstone.Lu2009Consistency,BJNP,Ma2013,CaiMaWu13}.

While standard (nonrestricted) PCA can be efficiently solved by
computing the eigenvectors of a symmetric matrix, sparse PCA is a
difficult combinatorial problem, and in fact solving ${\mathcal
{L}}_0(\Sigma)$ is
NP-hard.\setcounter{footnote}{3}\footnote{This claim follows from \cite{Moghaddam06generalizedspectral}
and \cite{Natarajan1995}, but can also proved by a direct reduction
from the
$k$-clique problem in a $p$-vertex graph,
and considering $\Sigma=A+pI$ where $A$ is the graph's adjacency matrix.
}
Nevertheless, various
computationally efficient approaches were developed to deal with the
problem. These include greedy or nonconvex optimization procedures
\cite{ZhangZhaSimon,JolliffeTrendafilovUddin03}, methods based on
$\ell
_1$-regularization \cite
{AspremontEtAlSDP07,ZouHastieTibshirani06,Moghaddam06spectralbounds,WittenTibshirani:2009},
regularized singular-value-decomposition \cite{Shen:2008}, an augmented
Lagrangian method \cite{ZhaosongZhang:2012}, a simple diagonal
thresholding (DT) algorithm \cite{Johnstone.Lu2009Consistency}, and
sophisticated semidefinite programming (SDP) methods \cite
{AspremontBannerjeeGhaoui07}. The latter approach, and in particular
its ability to recover an $\ell_0$-sparse PC, are the focus of the
current paper.

\textit{SDP-based algorithm}.
We study the following concrete SDP relaxation of \eqref{eq:LP},
which was suggested by d'Aspremont et al. \cite{AspremontEtAlSDP07}:
%
\begin{equation}
\label{eq:SDP} \operatorname{argmax} \bigl\{ \langle{\hat{\Sigma}}, X \rangle
\dvtx X \in\mathcal{S}^p_+, \operatorname{tr}(X)=1, \|X
\|_S \leq k\bigr\},
\end{equation}
where for two matrices $X,Y\in{\mathbb{R}}^{p\times p}$ we denote by
$\langle{X,Y}\rangle= \sum_{i,j} X_{ij}Y_{ij} = \operatorname
{tr}(X^TY)$ their Frobenius
inner-product, $\|X\|_S=\sum_{i,j}|X_{ij}|$ is the ``absolute-sum norm,''
and $\mathcal{S}^p_+ = \{X \in{\mathbb{R}}^{p \times p}\dvtx
X=X^T, X
\succeq0 \}$
is the cone of symmetric positive semidefinite (PSD) matrices.
As SDP \eqref{eq:SDP} returns a symmetric matrix rather than a vector,
d'Aspremont et al. \cite{AspremontEtAlSDP07} suggested to output its
leading eigenvector
as an estimate for the first sparse-PC. This algorithm is summarized as follows.


\begin{algorithm}[t]\label{alg:SDP}
\SetKwData{Left}{left}\SetKwData{This}{this}\SetKwData{Up}{up}
\SetKwFunction{Union}{Union}\SetKwFunction{FindCompress}{FindCompress}
\SetKwInOut{Input}{input}\SetKwInOut{Output}{output}

\Input{(mean-centered) vectors $\mathbf{x}_1,\ldots,\mathbf{x}_n
\in{\mathbb{R}}^p$, sparsity
parameter $k$}
\Output{vector $\hat\mathbf{z}\in{\mathbb{R}}^p$}
\caption{SDP-estimator}
\BlankLine
let ${\hat{\Sigma}}= \frac{1}{n}\sum_{i=1}^n \mathbf{x}_i \mathbf
{x}_i^T$ \\
compute a solution $X \in{\mathbb{R}}^{p \times p}$ of SDP \eqref
{eq:SDP} \\
let $\hat\mathbf{z}$ be the leading (unit-length) eigenvector of $X$
\end{algorithm}

\textit{Single-spike input model}.
We examine Algorithm \ref{alg:SDP} under
the single-spike multivariate Gaussian model introduced in \cite{Johnstone01},
where the samples $\mathbf{x}_i$ are of the form
%
\begin{equation}
\label{eq:defOfDist} {\mathbf{x}}_i = \sqrt{\beta}u_i {
\mathbf{z}} + \bolds\xi_i,\qquad i=1,\ldots,n.\vadjust{\goodbreak}
\end{equation}
Here, the parameter $\beta>0$ is the \emph{signal strength},
$\mathbf{z}\in{\mathbb{R}}^p$ is the \emph{planted spike} assumed
to be a $k$-sparse
unit-length vector,
$\bolds\xi_i \in{\mathbb{R}}^p$ is a noise vector whose entries are
all i.i.d.
 ${{N}}
(0,1)$ and $u_i\sim{{N}}(0,1)$. Furthermore, all the $u_i$'s and
$\bolds
\xi
_i$'s are independent of each other. The corresponding population
covariance matrix is
%
\begin{equation}
\Sigma= \beta\mathbf{z}\mathbf{z}^T + {\mathbf{I}}_p,
\label{eq:Sigma_rank_one}
\end{equation}
and its largest eigenvalue is $1+\beta$, with associated eigenvector
$\mathbf{z}$.
We consider throughout the scenario $(n,p,k)\to\infty$, and mention
additional assumptions (e.g., $\beta$ is fixed or $p/n$ tends to $c>0$)
as needed.

\textit{Information versus computational limits}.
Amini and Wainwright \cite{AminiWain09} studied this single-spike
input model, under the
additional assumption that the nonzero entries of $\mathbf{z}$ are
exactly of
the form $\pm1/\sqrt{k}$, which represents the hardest type of
$k$-sparse vectors.
They proved that up to sparsity level $k =O(\kappa_{n,p})$ where
$\kappa
_{n,p} = \sqrt{n/\log p}$,
Algorithm \ref{alg:SDP} outputs a vector $\hat\mathbf{z}$
whose support coincides with that of $\mathbf{z}$;\footnote{For technical
reasons, their proof requires the additional condition
$k = O(\log p)$, which they conjecture can be removed.}
they further showed, using a simple second moment calculation,
that up to the same order of sparsity level $k = O(\kappa_{n,p})$,
the diagonal thresholding algorithm \cite{Johnstone.Lu2009Consistency}
also recovers the support of $\mathbf{z}$
and fails whenever $k/\kappa_{n,p} \to\infty$.
In contrast, Amini and Wainwright \cite{AminiWain09} showed that for
$k=\Omega(\kappa
_{n,p}^2)$,%
\footnote{
We write $f=\Omega(g)$ if $f(n) \ge C g(n)$ for some absolute positive
constant $C$ and all sufficiently large $n$. Similarly, $f=\Theta(g)$
means $C_1 g(n) \leq f (n) \leq C_2 g(n)$.
}
\emph{every} method [including exhaustive search over all
$p\choose k$ subsets of size $k$] will err with probability at least
$1/2$. In fact, even the simpler task of \emph{detecting} the presence
of a spike is not possible for this range of parameters, as recently
proved in \cite{RigolletClique,RigolletCliqueCOLT}. For further results
including minimax rates, under more general sparsity models, see \cite
{LeiVu2013,CaiMaWu13,Samworth14}.

The following question thus remained open:
Does Algorithm \ref{alg:SDP}, which is more sophisticated
and computationally heavy, outperform the simple DT algorithm?
Specifically, are there intermediate sparsity levels $\kappa_{n,p} <
k<\kappa_{n,p}^2 $
(such that $k/\kappa_{n,p} \to\infty$, and ignoring multiplicative constants)
for which $\hat\mathbf{z}$ still approximates $\mathbf{z}$ in some
useful sense?
While not answering this question,
Amini and Wainwright proved that for sparsity level up to $k=O(\kappa
_{n,p}^2)$,
\emph{if} the solution to \eqref{eq:SDP} remains rank one,
\emph{then} the support of $\hat\mathbf{z}$ coincides with that of
$\mathbf{z}$.
They then suggested that for this $\ell_0$-sparse PCA problem, the
information and computational limits coincide, both are equal $\Theta
(\kappa_{n,p}^2)$,
and Algorithm \ref{alg:SDP} is optimal.
In their words, ``under the rank-one condition, the SDP is in fact
statistically optimal, that is, it requires only the necessary number
of samples (up to a constant factor) to succeed'' \cite{AminiWain09}, page
2880.

Our results, formally stated below,
prove that unfortunately this is \emph{not} the case---in fact,
when $k$ slightly exceeds $\kappa_{n,p}$, namely
$k = \Omega(\kappa_{n,p} \sqrt{\log p)}=\Omega(\sqrt{n})$,
the solution $X$ of SDP \eqref{eq:SDP} \emph{does not} have rank one
and is not close to $\mathbf{z}\mathbf{z}^T$.
Furthermore, if $X$ has a low rank,
then the output $\hat\mathbf{z}$ of Algorithm \ref{alg:SDP}
is at best weakly correlated with $\mathbf{z}$.
In Section~\ref{sec:Experimental} we present empirical simulation
results showing that indeed Algorithm \ref{alg:SDP} and DT perform similarly.

Given that the SDP algorithm does not seem to significantly improve
over DT under the single spike model, the following question arises:
Is there a simple algorithm which outperforms both? Motivated by the
work of Bickel and Levina \cite{BickelLevina06}, we suggest a
light-weight greedy algorithm called Covariance Thresholding (CT),
which can be seen as a generalization of Diagonal Thresholding. We
provide experimental results suggesting that CT is consistent for $k =
O(\sqrt{n})$; see Section~\ref{sec:Experimental} for details.
Recently, following our work, Deshpande and Montanari \cite
{MontanariCovar} rigorously proved that a variant of our CT algorithm
indeed asymptotically recovers the support of $\mathbf{z}$ up to these sparsity
levels. Finally, we note that despite our results, there are other
settings, such as estimating sparse eigenvectors of correlation
matrices, where SDP-based methods are provably better than diagonal
thresholding, possibly even achieving the relevant minimax rates \cite
{LeiVu14,WangLuLiu14}.

\subsection{Our results}\label{sec:OurResults}
We consider the single-spike model defined in \eqref{eq:defOfDist} in
high-dimensional settings whereby $(n,p,k)\to\infty$ and $p/n^\alpha
\to
c$ for positive constants $c,\alpha\ge1$. We further assume that the
$k$-sparse vector $\mathbf{z}$ has $k$ nonzero entries of the
form $\pm
1/\sqrt{k}$. In what follows, we denote by $\operatorname
{supp}(\mathbf{x})$ the set $\{
i\dvtx\mathbf{x}
_i \ne0\}$.
In the analysis, we assume without loss of generality that
the nonzero coordinates of the spike $\mathbf{z}$ are exactly its
first $k$
coordinates, that is, $\operatorname{supp}(\mathbf{z})=\{1,2,\ldots
,k\}$.

For the case $\alpha=1$, that is, $p/n\to c$, we focus on weak signal
strengths $\beta\leq\sqrt{\frac{p} n}$, whereas when $\alpha>1$, the
signal strength may grow to infinity provided it still satisfies $\beta
\leq\sqrt{\frac{p} n}$; see assumption (b) below. The reason is that
when $\alpha=1$ and $\beta>\sqrt{\frac{p} n}$, as the next theorem
shows, recovering the support of $\mathbf{z}$ is computationally
easy, almost
up to the information limit. As before, we let $\kappa_{n,p}=\sqrt
{n/\log p}$.

\begin{theorem}[(Strong signal)] \label{thm:largesignal}
Fix $c > 1$ and $\beta> \sqrt{c}$,
and let $(n,p,k)\to\infty$ such that $p/n\to c$ and $k/\kappa_{n,p}^2
\to0$.
Let $\hat\mathbf{w}_1$ be the leading eigenvector of ${\hat{\Sigma
}}$, and denote by
$\operatorname{supp}_k(\hat\mathbf{w}_1)$ its $k$ largest entries
in absolute value. Then
$\operatorname{supp}_k(\hat\mathbf{w}_1)=\operatorname
{supp}(\mathbf{z})$ with probability tending to one as
$(n,p,k)
\to\infty$.
\end{theorem}

Our next results, stated in the three theorems below,
refer to the following assumptions:
\begin{longlist}[(a)]
\item[(a)] 
Fix positive $c,\alpha\ge1$, and let $(n,p,k)\to\infty$ such that
$p/n^\alpha\to c$.
\item[(b)] 
The signal strength, either fixed or growing with $n,p$, satisfies
$\beta\le\sqrt{\frac{p}{n}}$.
\item[(c)] 
The sparsity level $k$ satisfies $k \ge2p/\sqrt{n}$, and $k/p \to0$.
\end{longlist}

We next analyze the quality of the output $\hat\mathbf{z}$ of
Algorithm \ref
{alg:SDP},
as measured by its cosine-similarity to the planted spike $\mathbf{z}$.

\begin{theorem}[(Cosine similarity)] \label{thm:main1}
Assume \textup{(a)--(c)}.
Then there exists ${\varepsilon}={\varepsilon}(n)\to0$,
such that if $X$ is a solution of SDP \eqref{eq:SDP}, and $\lambda_1$
is its largest eigenvalue,
then with probability tending to one as $(n,p,k)\to\infty$,
the output $\hat\mathbf{z}$ of Algorithm~\ref{alg:SDP} satisfies
%
\begin{equation}
\label{eq:DOT_PRODUCT}
\bigl\vert\langle{\hat\mathbf{z},\mathbf{z}}\rangle\bigr\vert^2
\le\frac{23}{\lambda_1}\sqrt{\frac{n}{p}} (1+\sqrt{\beta } ) +
\frac{{\varepsilon}}{\lambda_1}.
\end{equation}
\end{theorem}

The following corollary of Theorem~\ref{thm:main1} shows that the SDP
solution is far from $\mathbf{z}\mathbf{z}^T$. For a matrix $A$ we
denote its
spectral norm by $\|A\|=\sqrt{\lambda_{\max}(AA^T)}$.

\begin{corollary}\label{cor:NoCloseSolutions} Assume \textup{(a)--(c)}, and further that $p \ge150^4n$. Let $X$ be a solution of
SDP \eqref{eq:SDP}. Then $\| X -\mathbf{z}\mathbf{z}^T\| \ge\frac{1}{3}$ with
probability tending to one as $(n,p,k)\to\infty$.
\end{corollary}

\begin{pf}
Assume for contradiction that the matrix $Y=X-\mathbf{z}\mathbf{z}^T$
has a small
spectral norm $\eta_1 = \|Y\| < 1/3$.
Using Weyl's inequality \cite{STEWART90}, $\|X\|\geq\|\mathbf
{z}\mathbf{z}^T\| - \|
Y\|
$. Since $\|\mathbf{z}\|_2=1$, the largest eigenvalue of $X$ is thus lower
bounded by $\lambda_1 \ge1-1/3=2/3$. Let $\hat\mathbf{z}$ be a (unit-length)
eigenvector of $X$ corresponding to this largest eigenvalue $\lambda
_1$. Recalling the variational definition of
the largest eigenvector of a matrix, we obtain
%
\begin{equation}
\label{cor:eq1} \tfrac{2}{3} \le\lambda_1= \hat
\mathbf{z}^T X \hat\mathbf{z}= \hat\mathbf{z}^T \bigl(Y+
\mathbf{z}\mathbf{z}^T\bigr) \hat\mathbf{z}= \hat\mathbf{z}^T
Y \hat\mathbf{z}+ \hat\mathbf {z}^T \mathbf{z}\mathbf{z}^T
\hat\mathbf{z}.
\end{equation}
Using our assumption, $\hat\mathbf{z}^T Y \hat\mathbf{z}\le\|Y\| =
\eta_1 \le1/3$.
By Theorem~\ref{thm:main1}
%
\begin{equation}
\label{cor:eq2} \hat\mathbf{z}^T \mathbf{z}\mathbf{z}^T
\hat\mathbf{z}= \bigl\vert \langle{\hat\mathbf{z},\mathbf{z}}\rangle\bigr\vert^2
\le \frac
{23}{\lambda_1}\sqrt{\frac{n}{p}} (1+\sqrt{\beta} ) +
\frac{{\varepsilon}}{\lambda_1}.
\end{equation}
Plugging $p/n = 150^4$, $\beta\le\sqrt{p/n}=150^2$ and $\lambda_1
\ge2/3$
into equation (\ref{cor:eq2}) gives that its right-hand side is at most
$0.2315+\frac{3}2{\varepsilon}$.
Since by Theorem~\ref{thm:main1}, ${\varepsilon}={\varepsilon}(n)\to
0$ as $n\to
\infty$,
\eqref{cor:eq2} is strictly smaller than $1/3$ for a sufficiently
large $n$.
Combining \eqref{cor:eq1} and \eqref{cor:eq2} we arrive at the
following contradictory set of inequalities:
\[
\tfrac{2}{3} \le\lambda_1=\hat\mathbf{z}^T Y \hat
\mathbf{z}+ \hat \mathbf{z}^T \mathbf{z}^T\mathbf{z} \hat
\mathbf{z}< \tfrac{1}{3}+\tfrac{1}{3}=\tfrac{2}{3}.
\]
\upqed\end{pf}

Note that the constant 23 appearing in equation (\ref{eq:DOT_PRODUCT}),
and consequently the factor $150^{4}$ in the corollary, are not
necessarily optimal. Both may be further reduced at the expense of more
involved proofs.

Further note that if $\lambda_1(X)$ is bounded away from zero as $
p,n\to\infty$, then for $\alpha>1$, equation (\ref{eq:DOT_PRODUCT})
implies that $\langle{\hat\mathbf{z},\mathbf{z}}\rangle\to0$.
Namely, in this case the
output of Algorithm \ref{alg:SDP} is nearly orthogonal to $\mathbf
{z}$. Such an
empirical behavior of $\lambda_1$ was observed in our experimental
results; see Figure~\ref{fig:Alg1}.

We prove Theorem~\ref{thm:main1} using the next result, which itself
may be of interest as it bounds the value of SDP \eqref{eq:SDP}. Recall
that the SDP solution is highly nonlinear in its inputs, and therefore
no closed-form explicit expression is known for the solution $X$ or the
SDP value $\langle{\hat{\Sigma}}, X \rangle$.

\begin{theorem}[(SDP value)] \label{thm:main2}
Assume \textup{(a)--(c)}.
Then there exists $\zeta=\zeta(n) \to0$
such that with probability tending to one as $(n,p,k)\to\infty$,
every solution $X$ of SDP \eqref{eq:SDP} satisfies
%
\begin{equation}
(1-\zeta) \biggl(1+\frac{p}{n} \biggr) \leq\langle{\hat{\Sigma}}, X
\rangle \leq(1+\zeta) \biggl(1+\sqrt {\frac{p}{n}}+\sqrt{\beta}
\biggr)^2. \label{eq:BOUNDS_SDP}
\end{equation}
\end{theorem}

For $\alpha> 1$, the ratio between the upper and lower bounds
in \eqref
{eq:BOUNDS_SDP} is at most $1+O (\zeta+\sqrt{\frac
{n}{p}}(1+\sqrt
{\beta}) )$ and tends to one as $p,n\to\infty$.

For the important regime $\alpha=1$, we can use Theorem~\ref
{thm:main2} to sharpen our conclusion from Theorem~\ref{thm:main1} and
show that
with probability tending to one, not only $X \ne\mathbf{z}\mathbf
{z}^T$, but $X$ is
not even rank one.
We arrive at this conclusion by combining Theorem~\ref{thm:main2}
with the next theorem.

\begin{theorem} \label{thm:main3}
Assume \textup{(a)--(c)}, and in addition $\alpha= 1$, $c>20$ and
$k/  (p/\log^2 p ) \to0$. Then with probability tending to one
as $(n,p,k)\to\infty$, every rank-one matrix $Y=\mathbf{y}\mathbf
{y}^T$ that is feasible
for SDP \eqref{eq:SDP} satisfies
%
\begin{equation}
\langle{\hat{\Sigma}}, Y \rangle\le\frac{8}{9}\cdot
\frac{p}{n}. \label{eq:BOUND_RANK_ONE}
\end{equation}
\end{theorem}

To see that the solution $X$ of SDP \eqref{eq:SDP} is indeed not rank one,
we compare the upper bound in \eqref{eq:BOUND_RANK_ONE} with the
(larger) lower bound in \eqref{eq:BOUNDS_SDP}, namely,
$ \langle{\hat{\Sigma}}, Y \rangle
\le\frac{8}{9}\cdot\frac{p}{n}
< (1-\zeta) (1+\frac{p}{n} )
\leq\langle{\hat{\Sigma}}, X \rangle$.\footnote{We remark that another lower bound $\langle{\hat{\Sigma
}}, X \rangle \ge1+\beta$
was proved in \cite{RigolletClique}, Proposition~6.1,
in a setting similar to Theorem~\ref{thm:main2},
but we cannot use it to derive $ \langle{\hat{\Sigma}}, Y \rangle <
\langle{\hat{\Sigma}}, X \rangle$
because $\frac{8}{9} \frac{p}{n}$ could be larger than $1+\beta$.}

In conclusion, Theorems \ref{thm:main1}--\ref{thm:main3} suggest that
the standard SDP-based approach (provided by Algorithm \ref{alg:SDP})
is not significantly more effective than
the simpler, light-weight diagonal thresholding.
In particular, for weak signal strengths, Algorithm \ref{alg:SDP} does
not yield a rank-one solution and hence cannot provably solve sparse
PCA up to the information limit, as previously suspected.
Our conclusion is in line with a recent, independently obtained
result of Berthet and Rigollet \cite{RigolletCliqueCOLT}, which
asserts that
the existence of a polynomial-time computable statistic for reliably
\emph{detecting} the presence of a single spike of $\ell_0$-sparsity
$k$ for $k/\sqrt{n} \to\infty$,
implies a polynomial-time algorithm for reliably detecting the presence
of a planted clique of size $k'$, for $k'/\sqrt{n}\to0$, in an
otherwise random graph $G(n,1/2)$.
The latter problem, known as the \textit{hidden clique problem} in the
computer science literature, is believed to be a computationally hard
task, and polynomial-time algorithms known to date can only find a
planted clique whose size $k'$ is at least of order $\sqrt{n}$ \cite
{AKSClique,FK00,FeigeClique,Ames:2011,MontanariClique,PeresClique}.
Furthermore, Wang et al. \cite{Samworth14} showed that under the
hidden clique hardness assumption, in certain sparsity regimes no
randomized polynomial time algorithms can estimate the leading spiked
eigenvector with optimal rate.

Our result differs from \cite{RigolletCliqueCOLT} in several respects.
First, our results are unconditional; that is, Theorems \ref
{thm:main1}--\ref{thm:main3} are not based on any computational
hardness assumptions,
and thus remain valid even if future developments will yield
a polynomial-time algorithm for finding a hidden clique of size $n^{0.49}$.
Second, our focus is on estimation and not on detection, which in
general are different problems.

\begin{figure}

\includegraphics{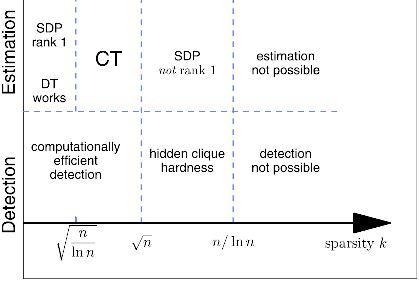}

\caption{State of the art for detection and estimation in the
single-spike model under various regimes of
$l_0$-sparsity (assuming $n \approx p$ and omitting constant
factors).}\label{fig:Ranges}
\end{figure}

We summarize in Figure~\ref{fig:Ranges}
the picture emerging from the results of Amini and Wainwright \cite
{AminiWain09}, Berthet and Rigollet
 \cite{RigolletCliqueCOLT}, Deshpande and Montanari \cite
{MontanariCovar} and our work.
Based on these results and the fact that even a sophisticated SDP-based
algorithm fails to estimate $\mathbf{z}$ for $k\ge\sqrt{n}$, we
conclude with
the following conjecture.

\begin{conjecture}
In the single-spike model with $p/n \to c$,
fixed signal strength $\beta\leq\sqrt{p/n}$
and $\ell_{0}$-sparsity $k=n^{0.5+{\varepsilon}}$ for fixed
${\varepsilon}> 0$,
no polynomial-time algorithm can recover the support of $\mathbf{z}$
with probability tending to one as $(n,p,k)\to\infty$.
\end{conjecture}

\textit{Organization}.
In Section~\ref{sec:CovarThresh} we describe our covariance
thresholding algorithm, followed by experimental results in
Section~\ref{sec:Experimental}. In Section~\ref
{sec:ProofThmLargeBeta} we give a
short proof of Theorem~\ref{thm:largesignal}. In Section~\ref
{sec:prelims} we assert preliminary facts that will be later used in
the proofs of Theorem~\ref{thm:main1} in Section~\ref{sec:ProfOfThm1},
Theorem~\ref{thm:main2} in Section~\ref{sec:ProfOfThm2} and Theorem~\ref
{thm:main3} in Section~\ref{sec:ProfOfThm3}.

\section{Covariance thresholding algorithm}\label{sec:CovarThresh}
Motivated by the work of Bickel and Levina \cite{BickelLevina06}, we
suggest Algorithm \ref{alg2} for the $\ell_0$-sparse PCA problem,
which we call \emph{covariance thresholding}, or CT for short.

\begin{algorithm}[t]
\SetKwData{Left}{left}\SetKwData{This}{this}\SetKwData{Up}{up}
\SetKwFunction{Union}{Union}\SetKwFunction{FindCompress}{FindCompress}
\SetKwInOut{Input}{input}\SetKwInOut{Output}{output}

\Input{vectors $\mathbf{x}_1,\ldots,\mathbf{x}_n \in{\mathbb{R}}^p$,
threshold $t$, sparsity
level $k$}
\Output{subset $S \subseteq[p]$ of cardinality $k$}
\caption{Covariance thresholding}\label{alg2}
\BlankLine
compute ${\hat{\Sigma}}= \frac{1}{n}\sum_{i=1}^n \mathbf{x}_i
\mathbf{x}_i^T$
\\
compute $T\in{\mathbb{R}}^{p\times p}$ by thresholding the entries of
${\hat
{\Sigma}}
$, namely,
\[
T_{ij} = \cases{ \hat\Sigma_{ij}, &\quad $\mbox{if $|\hat
\Sigma_{ij}| > t$}$; \vspace*{2pt}
\cr
0,&\quad $\mbox{otherwise}$}
\]
\\
let $\mathbf{w}\in{\mathbb{R}}^p$ be the leading eigenvector of $T$\\
let $S\subseteq[p]$ contain the $k$ coordinates of largest absolute
value in $\mathbf{w}$
\end{algorithm}

We present some intuition as to why we expect this algorithm to work.
From the definition of ${\hat{\Sigma}}$ in \eqref{eq:defOfDist},
it follows easily that the off-diagonal noise entries
have expected value zero and standard deviation $1/\sqrt{n}$,
while for signal entries the expected value is $\pm\beta/k$ with s.d.
$C(\beta)/\sqrt{n}$. Consider, for example, a signal strength $\beta=
1$, sparsity $k\le\sqrt{n}/10$ (where 10 is rather arbitrary), and
choose $t=5/\sqrt{n}$. Then for a noise entry to survive thresholding,
it must deviate from its mean by $5$ s.d. and an analogous deviation
for a signal entry to be zeroed out. Both events happen with small
constant probability; hence most noise entries are zeroed and a
constant fraction of signal entries survive.
In fact, when $k=O(\sqrt{n/ \log p})$ one can easily show that CT,
similar to DT, recovers the support of $\mathbf{z}$.
Recently, Deshpande and Montanari \cite{MontanariCovar} proved that a
variant of our algorithm
is consistent up to sparsity levels $k=O(\sqrt{n})$.
Their proof method is not directly applicable to our algorithm,
but simulation results, detailed below, suggest that our algorithm is also
able to recover the correct support up to $k=O(\sqrt{n})$.
Hence, covariance thresholding is thus far the only algorithm, with
polynomial run-time, that
can \textit{provably} recover the support up to sparsity levels $
k=O(\sqrt{n})$.

%

\section{Simulation results}\label{sec:Experimental}

\subsection{Covariance thresholding versus diagonal thresholding}
\begin{figure}

\includegraphics{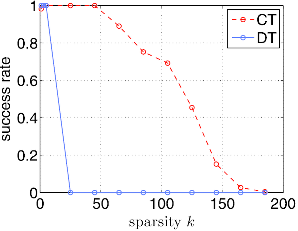}

\caption{Performance of DT vs. CT, $n=p=5000$.
$y$-axis is the success rate averaged over
500 runs, with signal strength $\beta=2$, and CT parameterized with
threshold $t=3/(2k)$.}\label{fig:CTvsDT}
\end{figure}

\begin{figure}[b]

\includegraphics{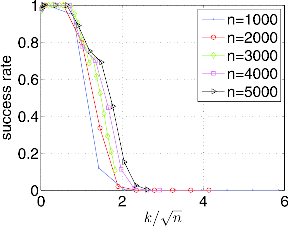}

\caption{Performance of CT in proportion to $k/\sqrt{n}$
(depicted for different $k$ and $n$).
$y$-axis is the success rate averaged over
500 runs, with signal strength $\beta=2$, and CT parameterized with
threshold $t=3/(2k)$.}\label{fig:CTScaled}
\end{figure}

We compare a few algorithms under the following setup. We generate $
n$ i.i.d. samples $\mathbf{x}_i$ from the single-spike model \eqref
{eq:defOfDist} with a spike $\mathbf{z}$ of the form $\mathbf
{z}= (\frac
{1}{\sqrt
{k}},\frac{1}{\sqrt{k}},\ldots,\break \frac{1}{\sqrt{k}},0,0,\ldots
, 0 )$.
We assume the sparsity level $k$ is a priori known, and say that an
execution of an algorithm is \emph{successful} if it returns the
support of $\mathbf{z}$ exactly, that is, if the output is the set $\{
1,\ldots
,k\}$. The \emph{success rate} of an algorithm in $M$ independent
executions is the number of times it is successful divided by $M$. In
each experiment we fix $n=p$ and for various values of $k$ we measure
the success rate averaged over $M=500$ independent executions.
Figure~\ref{fig:CTvsDT} compares the performance of our CT algorithm to
DT. It is evident from this figure that in our setting, CT outperforms
DT. Figure~\ref{fig:CTScaled} shows the success rate of CT as a
function of the sparsity level $k$ scaled by $\sqrt{n}$, plotted for
five different values of $n$. These results reinforce our prediction
that CT works up to sparsity levels proportional to $\sqrt{n}$ (perhaps
even slightly more).

\subsection{SDP (Algorithm \texorpdfstring{\protect\ref{alg:SDP}}{1}) versus diagonal
thresholding}
\begin{figure}

\includegraphics{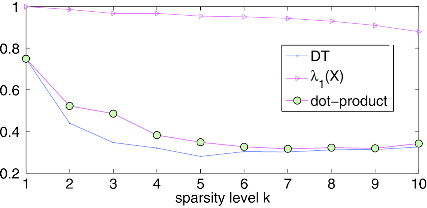}

\caption{Comparison of DT and SDP (Algorithm \protect\ref{alg:SDP})
for $n=p=50, \beta= 0.8$, averaged over 100 runs.
The blue dots represent the success rate of DT. The filled circles are
the average of the dot product
$|\langle\hat\mathbf{z},\mathbf{z}\rangle|$ of the SDP leading
eigenvector, whereas
the triangles are the
largest eigenvalue of the SDP solution, $\lambda_1(X)$.}
\label{fig:Alg1}
\end{figure}

We run Algorithm \ref{alg:SDP} with parameters $n=p=50$ and $\beta
=0.8$, averaging over $M=100$ runs. We solve the SDP in line 2 of
Algorithm \ref{alg:SDP} using SeDuMi 1.2.1 \cite{sedumi}.
Figure~\ref{fig:Alg1} plots the dot-product (in absolute value) between
$\hat\mathbf{z}$, the output of Algorithm \ref{alg:SDP} and the
planted spike
$\mathbf{z}$.
As expected, the dot-product gets smaller as the sparsity $k$ increases.
For comparison, the figure plots also the recovery rate of DT,
which also deteriorates as $k$ increases.
The figure also shows the largest eigenvalue of the SDP solution $X$;
we remark that this value is rather close to one,
even when the output of Algorithm~\ref{alg:SDP} is far from $\mathbf{z}$,
and is certainly bounded away from $0$, as assumed in the discussion
following Theorem~\ref{thm:main1}.

%

\section{Proof of Theorem \texorpdfstring{\protect\ref{thm:largesignal}}{1.1} (Strong signal)}
\label{sec:ProofThmLargeBeta}

Let $\hat\mathbf{w}_1$ be the leading eigenvector of ${\hat{\Sigma
}}$, and write it
as a linear combination of the spike $\mathbf{z}$ and some unit vector
$\mathbf{a}
\perp\mathbf{z}$,
namely, $\hat\mathbf{w}_1 = g\mathbf{z}+ \sqrt{1-g^2}\mathbf{a}$.
We may assume $g\in[0,1]$ by negating $\hat\mathbf{w}_1$, if necessary.
According to \cite{DPaul07}, Theorem~4,
for our setting of $\beta>\sqrt{c}$,
%
\begin{equation}
\label{eq:LargeBeta} g=g(\beta) \mathop{\rightarrow}^{\mathrm{a.s.}} \sqrt{\bigl(\beta
^2-c\bigr)/\bigl(\beta ^2+\beta c\bigr)} \qquad\mbox{as $n\to
\infty$}.
\end{equation}

Furthermore, according to Debashis (\cite{DPaul07}, Theorem~6), the
vector $\mathbf{a}\in{\mathbb{R}}^p$ is distributed uniformly on the unit
sphere of
dimension $p-1$ of vectors in $\mathbb{R}^p$ orthogonal to $\mathbf
{z}$.
Using this fact
we prove below the following property of the entries of~$\mathbf{a}$.

\begin{lemma}\label{lem:sphere} With probability tending to one,
all entries of $\mathbf{a}$ are bounded in absolute value by
$h\sqrt{\frac{\log p}{p}}$ for a suitable constant $h>0$.
\end{lemma}

Lemma~\ref{lem:sphere} implies that with probability tending to one,
for all $i \in[1,k]$ we have
$|(\hat\mathbf{w}_1)_i| \ge\frac{g}{\sqrt k}-\sqrt{1-g^2} \cdot
h\sqrt
{\frac
{\log p}{p}}$,
and for all $i \in[k+1,p]$ we have
$|(\hat\mathbf{w}_1)_i| \le\sqrt{1-g^2} \cdot h\sqrt{\frac{\log p}{p}}$.
To correctly identify the support of $\mathbf{z}$, it suffices to
require a gap
between signal and nonsignal coordinates, namely,
\[
\frac{g}{\sqrt{k}} > 2h \sqrt{1-g^2}\sqrt{\frac{\log p}{p}}.
\]
Solving for $k$ and using \eqref{eq:LargeBeta},
this inequality holds whenever $k<h' p/\log p$ for suitable
$h'=h'(\beta)>0$,
which in turn holds with probability tending to one,
because our assumption $k/\kappa_{n,p}^2 = k/(n/\log p) \to0$
implies $k/(p/\log p) \to0$.
This completes the proof of Theorem~\ref{thm:largesignal}.

\begin{pf*}{Proof of Lemma~\ref{lem:sphere}}
Let $\{\mathbf{s}_1,\ldots,\mathbf{s}_{p-1}\}$ be an orthonormal
basis for the
subspace of vectors in ${\mathbb{R}}^p$ orthogonal to $\mathbf{z}$. Since
$\mathbf{a}
=(a_1,\ldots
,a_p)$ is uniformly distributed in this subspace, it can be represented
as $\mathbf{a}= \frac{1}{\|\bolds\xi\|}\sum_{i=1}^{p-1} \xi
_i\mathbf{s}
_i$, where
$\bolds\xi=(\xi_1,\xi_2,\ldots,\xi_{p-1})$ is a vector of i.i.d.
standard Gaussians.

Fix a coordinate $i \in\{1,\ldots, p\}$, and write its corresponding
standard basis vector as $\mathbf{e}_i = \zeta_i \mathbf{z}+ \sqrt
{1-\zeta_i^2}
\tilde
\mathbf{e}_i$
for a unit vector $\tilde\mathbf{e}_i \perp\mathbf{z}$ and $\zeta
_i \in[-1,1]$.
Then $a_i = \mathbf{a}^T \mathbf{e}_i = \mathbf{a}^T (\zeta_i
\mathbf{z}+ \sqrt{1-\zeta_i^2}
\tilde
\mathbf{e}_i) = \sqrt{1-\zeta_i^2} \mathbf{a}^T \tilde\mathbf
{e}_i$, which implies
$|a_i| \le
|\mathbf{a}^T \tilde\mathbf{e}_i|$.
Since $\mathbf{a}$ and $\tilde\mathbf{e}_i$ are both unit vectors in
$\operatorname{span}\{\mathbf{s}
_1,\ldots,\mathbf{s}_{p-1}\}$, our task reduces to estimating the
inner-product between the uniformly distributed \emph{random} vector
$\mathbf{a}$ on the $(p-1)$-dimensional unit sphere and a \emph
{fixed} vector
$\tilde\mathbf{e}_i$ on the sphere.
Since $\mathbf{a}$ is random, we may replace $\tilde\mathbf{e}_i$
with another
fixed vector,
say $\mathbf{s}_{1}$. Namely, $ \mathbf{a}^T \tilde\mathbf{e}_i$
has the
same distribution
as $\mathbf{a}^T \mathbf{s}_1 = \xi_1/\|\bolds\xi\|$.
Standard tail bounds for the Gaussian and $\chi^2$ distributions (note
that $\Vert\bolds\xi\Vert^2 \sim\chi^2_{p-1}$) imply that $\frac
{|\xi
_1|}{\|
\bolds\xi\|}\le h\sqrt{\frac{\log p}{p}}$ with probability at least
$1-1/p^4$, for a suitable constant $h>0$.
The lemma follows by a union bound over all $p$ coordinates of $\mathbf{a}$.
\end{pf*}

\section{Preliminaries} \label{sec:prelims}
In this section we record a few standard results that will be used
later in the proofs.
The first is a large deviation result for a Chi-square random variable.

\begin{lemma}[(\cite{LAURENT-MASSART90})]\label{lem:XiSquare}
Let $X \sim\chi^2_{n}$. For all $x \geq0$,
\[
\operatorname{Pr}[X \geq n+2\sqrt{nx} + x] \leq e^{-x} \quad\mbox{and}\quad
\operatorname{Pr}[X \leq n-2\sqrt{nx}] \leq e^{-x}.
\]
\end{lemma}

The second lemma records a well-known argument about the inner-product
of two high-dimensional Gaussians.

\begin{lemma}\label{cl:CLT}
Let $\{x_i,y_i\}_{i=1}^n$ be standard i.i.d. Gaussian random variables.
Then $\sum_{i=1}^n x_i y_i$ is distributed like the product
of two independent random variables $\|\mathbf{x}\|\cdot\tilde y$,
where $\mathbf{x}
=(x_1,\ldots,x_n)$, $\Vert\mathbf{x}\Vert^2 \sim\chi^2_n$ and
$\tilde
y$ is a
standard Gaussian.
\end{lemma}

\begin{pf}
For every fixed realization of $\mathbf{x}$,
we have $x_i y_i \sim{{N}}(0,x_i^2)$ and by the independence of the $y_i$'s,
\[
\sum_{i=1}^n x_iy_i
\sim 
{{N}}\bigl(0,\Vert\mathbf{x}\Vert^2\bigr)=\Vert
\mathbf{x}\Vert\cdot {{N}}(0,1):=\Vert\mathbf{x}\Vert\cdot\tilde y.
\]
The lemma follows by observing that $\Vert\mathbf{x}\Vert^2 \sim
\chi^2_n$.
\end{pf}


The next proposition establishes an upper bound on $\lambda_{\max
}({\hat{\Sigma}})$, the maximal eigenvalue of the sample covariance matrix
${\hat{\Sigma}}$, in the single-spike model, in two regimes:
(i) $p/n^\alpha\to c$ for positive $c,\alpha\ge1$ and
(ii) $p/n \to0$.
The spectrum of the covariance matrix has been studied extensively in
the literature. Specifically, both Baik and Silverstein \cite{BaikSilverstein06},
Theorem~1.2 and Johnstone \cite{Johnstone01},
Theorem~1.1, provide the limiting behavior of $\lambda
_{\max}({\hat{\Sigma}})$ for $p/n\to c \ge1$ (i.e., $\alpha= 1$).
The regime of a fixed $p$ with $n\to\infty$ which implies $p/n
\to
0$ was analyzed in \cite{JoliffePCA}, Chapter~3, for example. Since we
could not locate a reference for the case $p/n^\alpha\to c$ and
$\alpha> 1$, or for $p/n \to0$ and $p$ not necessarily fixed, we
provide the following proposition.
The proof uses standard arguments and is given in
Section~\ref{sec:PrelimProofs}.

\begin{proposition}\label{prop:LargestEig} Let ${\hat{\Sigma}}$ be a $p
\times
p$ sample covariance matrix of $n$ samples in the $k$-sparse
single-spike model with signal strength $\beta> 0$, arbitrary $k$ and either:
\textup{(i)} $p/n \to0$ or
\textup{(ii)} $p/n^\alpha\to c$ for positive constants $c,\alpha\ge1$.
Then there exists an ${\varepsilon}= {\varepsilon}(n) \to0$ such
that with probability
tending to one as $n \to\infty$,
%
\begin{equation}
\label{eq:PropLargestEig} \lambda_{\max}({\hat{\Sigma}}) \leq(1+{\varepsilon})
\biggl(1+\sqrt {\frac{p}{n} }+\sqrt{\beta } \biggr)^2.
\end{equation}
\end{proposition}

\begin{corollary}\label{cor:TrSigHatzz} Let ${\hat{\Sigma}}$ be a $p
\times p$
sample covariance matrix of $n$ samples and a $k$-sparse spike $\mathbf{z}$
with signal strength $\beta> 0$. Further assume that $k/n \to0$.
Then there exists an ${\varepsilon}= {\varepsilon}(n) \to0$ such
that with probability
tending to one as $n \to\infty$, for every rank-one trace-one $p
\times p$ matrix $Y=\mathbf{y}\mathbf{y}^T$ with $\operatorname
{supp}(\mathbf{y})\subseteq\operatorname{supp}(\mathbf{z})$,
\[
\langle{\hat{\Sigma}}, Y \rangle \leq(1+{\varepsilon}) (1+\sqrt\beta
)^2.
\]
\end{corollary}

\begin{pf}
Consider $\sup_Y \langle{\hat{\Sigma}}, Y \rangle$ where $Y$
ranges over all matrices $Y$ as
stated above. For each such $Y=\mathbf{y}\mathbf{y}^T$, we have
$\|\mathbf{y}\|^2 = \sum_i \mathbf{y}_i^2 = \sum_i Y_{ii} =
\operatorname{tr}(Y) = 1$.
Let $\mathbf{y}_\mathbf{z}\in{\mathbb{R}}^k$ be the projection of
$\mathbf
{y}\in{\mathbb{R}}^p$ on the
coordinates of $\operatorname{supp}(\mathbf{z})$, then $\|\mathbf
{y}_\mathbf{z}\| = \|\mathbf{y}\|=1$.
Similarly, let
${\hat{\Sigma}}_\mathbf{z}$ be the $k \times k$ submatrix of ${\hat
{\Sigma}}$ corresponding
to $\operatorname{supp}(\mathbf{z})$, namely, restricting it to the
first $k$ rows and first
$k$ columns.
Observe that we can write
\[
\langle{\hat{\Sigma}}, Y \rangle=\operatorname{tr}\bigl({\hat{\Sigma }}
\mathbf{y}\mathbf{y}^T\bigr) = \mathbf{y}^T {\hat{\Sigma}}
\mathbf{y}= \mathbf{y}_\mathbf{z}^T {\hat{\Sigma}}
_\mathbf{z}\mathbf{y}_\mathbf{z} \leq\lambda_{\max}({\hat{
\Sigma}}_\mathbf{z}),
\]
hence $\sup_Y \langle{\hat{\Sigma}}, Y \rangle \leq\lambda_{\max
}({\hat{\Sigma}}_\mathbf{z})$.
Now the desired upper bound on $\lambda_{\max}({\hat{\Sigma
}}_\mathbf{z})$
follows using the fact $k/n \to0$ from Proposition~\ref{prop:LargestEig},
that is, plugging $p=k$ into \eqref{eq:PropLargestEig}.
\end{pf}

Our next proposition estimates $\operatorname{tr}({\hat{\Sigma}})$
and $\operatorname{tr}({\hat{\Sigma}}^2)$ for
the case $\beta= 0$ (no signal). These estimates were derived in \cite{Ledoit02},
Proposition~1, for example, but again only for $\alpha=
1$. For lack of reference we reprove it for $\alpha\ge1$
in Section~\ref{sec:PrelimProofs}.

\begin{proposition}\label{fact:ledoit}
Let ${\hat{\Sigma}}$ be a $p \times p$ sample covariance matrix of $n$
multivariate Gaussian observations whose population covariance matrix
is the identity. Assume that $(\log p)/n \to0$ as $n,p \to\infty$.
Then there exists an ${\varepsilon}={\varepsilon}(n) \to0$ such that
with probability
tending to one as $n\to\infty$,
\begin{eqnarray*}
(1-{\varepsilon})p\le\operatorname{tr}({\hat{\Sigma}})& \le& (1+{\varepsilon})p,
\\
(1-{\varepsilon})p \biggl(1+\frac{p}{n} \biggr) &\le&\operatorname {tr}
\bigl({\hat {\Sigma}}^2\bigr) \le(1+{\varepsilon} )p \biggl(1+
\frac{p}{n} \biggr).
\end{eqnarray*}
\end{proposition}

\section{Proof of Theorem \texorpdfstring{\protect\ref{thm:main1}}{1.2} (Cosine similarity)}
\label{sec:ProfOfThm1}

Let $X$ be a solution to SDP~\eqref{eq:SDP},
with eigenvalues $\lambda_1 \ge\lambda_2 \ge\cdots\ge\lambda
_p\ge0$
and a corresponding orthonormal set of eigenvectors
$\hat\mathbf{z}=\mathbf{v}_1,\ldots,\mathbf{v}_p\in{\mathbb{R}}^p$.
We can then write $X= \sum_{i=1}^p\lambda_i \mathbf{v}_i \mathbf{v}_i^T$,
and by linearity of the Frobenius inner-product,
$
\langle{\hat{\Sigma}}, X \rangle
= \sum_{i=1}^p\lambda_i \mathbf{v}_i^T {\hat{\Sigma}}\mathbf{v}_i
$.
Using the simple observations
$\mathbf{v}_i^T{\hat{\Sigma}}\mathbf{v}_i \leq\lambda_{\max
}({\hat{\Sigma}})$
(by the variational characterization of eigenvalues)
and $\sum_i \lambda_i = \operatorname{tr}(X)=1$, we get
%
\begin{equation}
\label{eq:PeeXopt} \langle{\hat{\Sigma}}, X \rangle = \sum
_{i=1}^p\lambda_i \mathbf{v}_i^T
{\hat{\Sigma}}\mathbf{v}_i 
 \leq
\lambda_1 \hat\mathbf{z}^T{\hat{\Sigma}}\hat\mathbf{z}+
(1-\lambda_1)\cdot \lambda _{\max}({\hat{\Sigma}}).
\end{equation}

Let us first provide a high-level description of the proof idea.
We can bound $\langle{\hat{\Sigma}}, X \rangle$ from below
(using Theorem~\ref{thm:main2}, which we prove in Section~\ref
{sec:ProfOfThm2}, and as mentioned earlier is used here)
and $\lambda_{\max}({\hat{\Sigma}})$ from above
(using Proposition~\ref{prop:LargestEig})
both by roughly $\frac{p}{n}$.
Now suppose $\lambda_1$ is not too small;
then on the right-hand side of \eqref{eq:PeeXopt}, a large contribution
must come from the first term $\lambda_1 \hat\mathbf{z}^T {\hat
{\Sigma}}\hat\mathbf{z}$.
But the quadratic form $\hat\mathbf{z}^T {\hat{\Sigma}}\hat\mathbf
{z}$ has small value
in the direction $\hat\mathbf{z}=\mathbf{z}$ (using Corollary~\ref
{cor:TrSigHatzz}),
and thus $\hat\mathbf{z}$ and $\mathbf{z}$ cannot be too close to
each other.

We now proceed to the detailed proof,
starting with a lower bound on $\lambda_1 \hat\mathbf{z}^T {\hat
{\Sigma}}\hat\mathbf{z}$.
Assume henceforth that the high-probability event asserted
by Theorem~\ref{thm:main2} indeed occurs;
namely, inequality \eqref{eq:BOUNDS_SDP} holds.
Similarly Corollary~\ref{cor:TrSigHatzz}
implies that inequality \eqref{eq:PropLargestEig} holds.
Plugging these two bounds into \eqref{eq:PeeXopt}
and using $\lambda_1 > 0$, we get
\begin{eqnarray*}
\lambda_1 \hat\mathbf{z}^T {\hat{\Sigma}}\hat\mathbf{z}
&\geq&(1-\zeta)\frac{p}{n} - (1-\lambda_1) (1+{\varepsilon})
\biggl(1+\sqrt{\frac
{p}{n}}+\sqrt{\beta } \biggr)^2
\\
&\geq&\frac{p}{n}(\lambda_1-\zeta-{\varepsilon}) -
(1+{\varepsilon}) \biggl(2\sqrt{\frac{p}{n}} + 1 +\sqrt{\beta } \biggr) (1+
\sqrt{\beta} ).
\end{eqnarray*}
Observe that $\sqrt{\frac{p}{n}} \ge\frac{1}{1+{\varepsilon}_1}$
for suitable ${\varepsilon}_1\to0$ and sufficiently large $n,p$.
In addition, assumption (b) yields that
$\beta\leq\sqrt{\frac{p}{n}} \leq(1+{\varepsilon}_1)\frac{p}{n}$.
For suitable ${\varepsilon}_2=O(\zeta+{\varepsilon}+{\varepsilon
}_1)$, we get
%
\begin{equation}
\label{eq:DotProdLower} \lambda_1 \hat\mathbf{z}^T {\hat{\Sigma}}
\hat\mathbf{z} 
\geq\frac{p}{n}(
\lambda_1-{\varepsilon}_2) - (4+{\varepsilon}_2)
\sqrt{\frac{p}{n}} (1+\sqrt{\beta } ).
\end{equation}

Next, we analyze the quadratic form $\hat\mathbf{z}^T {\hat{\Sigma
}}\hat\mathbf{z}$
in terms of $\gamma=\langle{\hat\mathbf{z},\mathbf{z}}\rangle\in[-1,1]$.
Write $\hat\mathbf{z}= \gamma\mathbf{z}+\sqrt{1-\gamma^2}\mathbf{s}$,
where $\mathbf{s}$ is a unit vector orthogonal to $\mathbf{z}$,
and recall that our goal is to upper bound $\gamma^2$.
Using Cauchy--Schwarz and the triangle inequality,
%
\begin{eqnarray}
\label{eq:DotProdUpper1} \hat\mathbf{z}^T {\hat{\Sigma}}\hat\mathbf{z} &\le&
\Vert\hat\mathbf{z}\Vert \cdot\Vert{\hat{\Sigma}}\hat \mathbf{z}\Vert = 1\cdot
\bigl\Vert{\hat{\Sigma}}\bigl(\gamma\mathbf{z}+\sqrt{1-\gamma ^2}\mathbf{s}
\bigr)\bigr \Vert
\nonumber
\\[-8pt]
\\[-8pt]
\nonumber
&\le&\vert\gamma\vert\cdot\|{\hat{\Sigma}}\mathbf{z}\|+\sqrt {1-
\gamma^2}\|{\hat{\Sigma}} \mathbf{s}\|.
\end{eqnarray}
Since ${\hat{\Sigma}}$ is PSD, it can be written as ${\hat{\Sigma
}}= B^TB$
for some matrix $B$ whose spectral norm
is $\Vert B\Vert=\Vert B^T\Vert=\sqrt{\lambda_{\max}({\hat{\Sigma}})}$.
Assume henceforth that the high-probability event asserted
by Corollary~\ref{cor:TrSigHatzz} indeed occurs, and we have
$\|B\mathbf{z}\|^2 = \mathbf{z}^T {\hat{\Sigma}}\mathbf{z}=
\langle{\hat{\Sigma}}, \mathbf{z}\mathbf{z}^T \rangle \le
(1+{\varepsilon}
_3)(1+\sqrt{\beta})^2$
for suitable ${\varepsilon}_3\to0$.
Using Proposition~\ref{prop:LargestEig} similarly yields
$\Vert B\Vert^2
= \lambda_{\max}({\hat{\Sigma}})
\leq(1+{\varepsilon}_4) ( 1+\sqrt{\frac{p}{n}} + \sqrt\beta
 )^2$
for suitable ${\varepsilon}_4\to0$.
Together, for suitable ${\varepsilon}_5=O({\varepsilon
}_3+{\varepsilon}_4+{\varepsilon}_1)$,
\[
\Vert{\hat{\Sigma}}\mathbf{z}\Vert 
\le\bigl\Vert B^T\bigr\Vert
\cdot\Vert B\mathbf{z}\Vert 
\le(3+{\varepsilon}_5) \sqrt{
\frac{p}{n}} ( 1+\sqrt{\beta} ),
\]
and similarly
\[
\Vert{\hat{\Sigma}}\mathbf{s}\Vert \leq\lambda_{\max}({\hat{\Sigma}})
\cdot\Vert\mathbf{s}\Vert 
 \leq\frac{p}{n}(1+{
\varepsilon}_4)+(4+{\varepsilon}_5)\sqrt {
\frac{p}{n}} ( 1+\sqrt {\beta } ).
\]
Plugging these into \eqref{eq:DotProdUpper1}
and using $\vert\gamma\vert\leq1$ and
$\sqrt{1-\gamma^2} \leq1-\frac{\gamma^2}{2}\le1$,
we have
%
\begin{equation}
\label{eq:DotProdUpper2} \hat\mathbf{z}^T {\hat{\Sigma}}\hat\mathbf{z} \le
\biggl( 1-\frac{\gamma^2}{2} + {\varepsilon}_4 \biggr)
\frac{p}{n} + (7+2{\varepsilon}_5) \sqrt{\frac{p}{n}} ( 1+
\sqrt{\beta} ).
\end{equation}

Now combining this upper bound \eqref{eq:DotProdUpper2}
with our lower bound \ref{eq:DotProdLower} (after dividing by $\lambda
_1$), gives
\begin{eqnarray*}
&&\biggl( 1-\frac{{\varepsilon}_2}{\lambda_1} \biggr) \frac{p}{n} - \frac{4+{\varepsilon}_2}{\lambda_1}
\sqrt{\frac{p}{n}} (1+\sqrt{\beta } )\\
&&\qquad \leq \biggl( 1-\frac{\gamma^2}{2} +{
\varepsilon}_4 \biggr) \frac{p}{n} + (7+2{\varepsilon}_5)
\sqrt{\frac{p}{n}} ( 1+\sqrt{\beta} ),
\end{eqnarray*}
and by further manipulation,
\[
\frac{\gamma^2}{2} - \frac{{\varepsilon}_2}{\lambda_1} - {\varepsilon}_4
\le\frac{11 + 2{\varepsilon}_5 + {\varepsilon}_2}{\lambda_1} \sqrt{\frac{n}p} ( 1+\sqrt{
\beta} ).
\]
For sufficiently large $n,p$,
this yields the bound on $\gamma^2 = \vert\langle{\hat\mathbf
{z},\mathbf{z}}\rangle\vert^2$
asserted in \eqref{eq:DOT_PRODUCT}, and completes the proof of
Theorem~\ref{thm:main1}.

\section{Proof of Theorem \texorpdfstring{\protect\ref{thm:main2}}{1.4} (SDP value)}
\label{sec:ProfOfThm2}
We start with the upper bound on $\langle{\hat{\Sigma}}, X \rangle$.
The idea is to drop the constraint $\|X\|_S \leq k$ from SDP \eqref{eq:SDP},
and show that the value of the resulting SDP, which can only be bigger,
is actually $\lambda_{\max}({\hat{\Sigma}})$,
and is thus bounded by Proposition~\ref{prop:LargestEig}.

Formally, let $X$ be a solution to SDP \eqref{eq:SDP},
and let us argue that (with probability~$1$)
\[
\langle{\hat{\Sigma}}, X \rangle \leq\sup \bigl\{\langle{\hat{\Sigma}}, Y
\rangle \dvtx Y\in \mathcal{S}^p_+, \operatorname {tr}(Y)=1 \bigr\} =
\lambda_{\max}({\hat{\Sigma}}).
\]
Indeed, the inequality holds because we have just relaxed SDP \eqref{eq:SDP}.
The equality holds by the following standard argument.
Writing $Y=\sum_i \mu_i \mathbf{y}_i \mathbf{y}_i^T$,
where $\{\mu_i\}_i$ are the eigenvectors of $Y$
and $\{\mathbf{y}_i\}_i$ is a corresponding orthonormal eigenbasis,
we have
%
\begin{eqnarray*}
\langle{\hat{\Sigma}}, Y \rangle 
&= &\sum_i
\mu_i \mathbf{y}_i^T {\hat{\Sigma}}
\mathbf{y}_i \leq\lambda_{\max}({\hat{\Sigma}})\cdot\sum
_i \mu_i = \lambda_{\max}({
\hat{\Sigma}})\cdot\operatorname{tr}(Y) \\
&=& \lambda_{\max}({\hat{\Sigma}}),
\end{eqnarray*}
and equality is achieved when maximizing over all relevant $Y$,
by taking $Y=\mathbf{y}_1 \mathbf{y}_1^T$ to be a rank-one matrix
where $\mathbf{y}_1$ is a leading eigenvector of ${\hat{\Sigma}}$.

To conclude the upper bound asserted in the theorem,
we combine the above with Proposition~\ref{prop:LargestEig},
and get that for a suitable ${\varepsilon}={\varepsilon}(n)\to0$
with probability tending to one as $n \to\infty$,
\[
\langle{\hat{\Sigma}}, X \rangle \leq\lambda_{\max}({\hat{\Sigma}})
\le(1+{\varepsilon}) \biggl(1+\sqrt{\frac{p}{n}}+\sqrt{\beta }
\biggr)^2.
\]

We turn to proving the lower bound on $\langle{\hat{\Sigma}}, X
\rangle$.
The idea is to consider a specific $X^*$
which is feasible (but not necessarily optimal) for SDP \eqref{eq:SDP},
and compute its objective value $\langle{\hat{\Sigma}}, X^* \rangle$.
Our $X^*$ is based on taking the nonsignal part of ${\hat{\Sigma}}$
[which is a $(p-k)\times(p-k)$ submatrix], padded with zeros elsewhere,
and ``forcing'' it to satisfy the constraints of SDP \eqref{eq:SDP}
by scaling it to be trace-one.

Formally, let $X^*=\tilde\Sigma/\operatorname{tr}(\tilde\Sigma)$,
where the matrix $\tilde\Sigma\in{\mathbb{R}}^{p \times p}$ is
given by
%
\begin{equation}
\label{eq:R} \tilde\Sigma_{ij} = \cases{ %
 0, &\quad
$\mbox{if $i \leq k$ or $j \leq k$};$
\vspace*{2pt}\cr
{\hat{\Sigma}}_{ij}, & \quad$\mbox{otherwise.}$ }
\end{equation}
We prove below that with probability tending to one,
the following inequalities hold for a suitable $\zeta=\zeta(n)\to0$:
%
\begin{eqnarray}
\bigl\Vert X^*\bigr\Vert_S &\le&\frac{2p}{\sqrt{n}} \leq k, \label{eq:TrSigHatX}
\\
\bigl\langle{\hat{\Sigma}}, X^* \bigr\rangle &\ge&(1-\zeta) \biggl(1+
\frac
{p}{n} \biggr). \label{eq:NormOneOfX}
\end{eqnarray}
Combining this with $X^*\in\mathcal{S}^p_+$ and $\operatorname{tr}(X^*)=1$,
which hold by construction,
will prove that with probability tending to one,
$X^*$ is feasible and has a high-objective value.

Before proceeding to prove \eqref{eq:NormOneOfX} and \eqref{eq:TrSigHatX},
we observe that the nonzeroed part of $\tilde\Sigma$ satisfies the
conditions of Proposition~\ref{fact:ledoit}, as it is a $(p-k)\times
(p-k)$ sample covariance matrix of $n$ multivariate Gaussian
observations whose population covariance matrix is the identity,
and furthermore $(\log(p-k))/n \to0$.
The concrete bounds that we get hold for a suitable ${\varepsilon
}={\varepsilon}(n)\to0$
and with probability tending to one as $n \to\infty$,
and roughly say that $\operatorname{tr}(\tilde\Sigma)\approx p-k$
and $\operatorname{tr}(\tilde\Sigma^2)\approx(p-k)(1+\frac{p-k}{n})$.

Let us now prove inequality \eqref{eq:TrSigHatX}.
First, using Cauchy--Schwarz,
\[
\bigl\Vert X^*\bigr\Vert_S = \frac{1}{\operatorname{tr}(\tilde\Sigma)} \sum
_{i,j>k} \vert \tilde{\Sigma}_{ij} \vert \leq
\frac{\sqrt{(p-k)^2\sum_{i,j>k}\tilde{\Sigma}_{ij}^2}}{
\operatorname{tr}(\tilde\Sigma)} = (p-k) \frac{\sqrt{\operatorname{tr}(\tilde{\Sigma
}^2)}}{\operatorname{tr}(\tilde\Sigma)}.
\]
By the above bounds from Proposition~\ref{fact:ledoit},
with probability tending to one,
\begin{eqnarray*}
\frac{\sqrt{\operatorname{tr}(\tilde{\Sigma}^2)}} {\operatorname
{tr}(\tilde\Sigma)} &\le&
\frac{\sqrt{ (1+{\varepsilon})(p-k)(1+{(p-k)}/{n}) }}{
(1-{\varepsilon})(p-k)} \le\frac{\sqrt{(1+{\varepsilon})
({1}/{(p-k)}+{1}/{n}) }}{
1-{\varepsilon}}\\
& \le&
\frac{1+{\varepsilon}}{1-{\varepsilon}} \sqrt{\frac{3}{n}} ,
\end{eqnarray*}
which together imply that
$
\Vert X^*\Vert_S
\le(p-k) \frac{1+{\varepsilon}}{1-{\varepsilon}} \sqrt{\frac{3}{n}}
\le\frac{2p}{\sqrt n}
\le k $.

We next prove inequality \eqref{eq:NormOneOfX}.
First, we expand
\[
\bigl\langle{\hat{\Sigma}}, X^* \bigr\rangle = \frac{1}{\operatorname{tr}(\tilde\Sigma)} \sum
_{i,j} \tilde\Sigma_{ij}{\hat{\Sigma}}_{ij}
= \frac{1}{\operatorname{tr}(\tilde\Sigma)} \sum_{i,j} \tilde
\Sigma_{ij}^2 = \frac{\operatorname{tr}(\tilde\Sigma^2)} {\operatorname
{tr}(\tilde\Sigma)}.
\]
By the above bounds from Proposition~\ref{fact:ledoit},
with probability tending to one,
\[
\frac{\operatorname{tr}(\tilde\Sigma^2)}{\operatorname{tr}(\tilde
\Sigma)} \ge\frac{(1-{\varepsilon})(p-k) (1+{p}/{n}) (1-{k}/{p}) }{
(1+{\varepsilon})(p-k)} \ge(1-\zeta) \biggl(1+\frac{p}{n}
\biggr),
\]
for a suitable $\zeta=\zeta(n)\to0$,
where we used here that $k/p \to0$ by assumption~(c).
Altogether, we conclude that
$\langle{\hat{\Sigma}}, X^* \rangle \geq(1-\zeta) (1+\frac
{p}{n} )$.

Having proved inequalities \eqref{eq:NormOneOfX} and \eqref{eq:TrSigHatX},
we conclude that with probability tending to one,
$X^*$ is feasible and has a high objective value,
which establishes a lower bound on the optimal SDP value $\langle{\hat
{\Sigma}}, X \rangle$,
and completes the proof of Theorem~\ref{thm:main2}.

\section{Proof of Theorem \texorpdfstring{\protect\ref{thm:main3}}{1.5} (SDP value)}
\label{sec:ProfOfThm3}

Let $F$ be the set of all vectors $\mathbf{y}$ whose corresponding
rank-one matrix
$Y=\mathbf{y}\mathbf{y}^T$ is feasible for SDP \eqref{eq:SDP}, formally,
\[
F = \bigl\{ \mathbf{y}\in{\mathbb{R}}^p \dvtx\| \mathbf{y}
\|_2 \leq1 \mbox{ and } \| \mathbf{y}\|_1 \leq \sqrt{k}
\bigr\}.
\]
We need to prove that with probability tending to one as $n \to\infty$,
every $Y=\mathbf{y}\mathbf{y}^T$ such that $\mathbf{y}\in F$
satisfies $\langle{\hat{\Sigma}}, Y \rangle \le\frac{8}{9}\frac{p}{n}$.
At a high level, $\langle{\hat{\Sigma}}, Y \rangle=\mathbf
{y}^T{\hat{\Sigma}}\mathbf{y}$ is
continuous,
and thus a standard approach is to discretize $F$ with an ${\varepsilon}$-net,
analyze every single point in $F$ separately
and apply a union bound argument.
The size of an ${\varepsilon}$-net for the unit $\ell_2$-ball in $p$
dimensions
is proportional to $(1/{\varepsilon})^p$.
On the other hand, our upper bound on the probability that a fixed
$Y=\mathbf{y}
\mathbf{y}^T$
violates $\langle{\hat{\Sigma}}, Y \rangle\le\frac{8}{9}\frac
{p}{n}$ is larger than ${\varepsilon}^{p}$; see
Lemma~\ref{lem:BounOnQuotient}.
Therefore, a naive discretization of $F$ fails,
and we need to reduce the size of the net by using the additional constraint
$\Vert\mathbf{y}\Vert_1 \leq\sqrt{k}$.
To this end, we approximate $F$ by a set $\hat F$,
whose definition uses an $\ell_0$-constraint;
the idea is that an $\ell_0$-bound is technically more convenient than
$\ell_1$.
We apply an ${\varepsilon}$-net argument to $\hat F$,
which indirectly yields a bound for all of $F$.
Specifically, we define
\[
\hat F = \bigl\{ \mathbf{y}\in{\mathbb{R}}^p \dvtx\Vert\mathbf{y}
\Vert _2 \leq1 \mbox{ and } \Vert\mathbf{y}\Vert_0 \leq40
\sqrt{pk} \bigr\}.
\]
To formalize the notion of one set approximating another one,
we define the \emph{$r$-neighborhood} of a set $A\subset{\mathbb
{R}}^p$ to be
$A_r = \{ \mathbf{y}\in{\mathbb{R}}^p \dvtx \exists\mathbf{y}'
\in A, \|
\mathbf{y}-\mathbf{y}'\| \leq r
\}$.

\begin{lemma}\label{lem:epsNet}
The sets $\hat F,F$ defined above satisfy $F \subseteq\hat F_{1/40}$.
\end{lemma}

\begin{pf}
Fix $\mathbf{y}\in F$, and let $I = \{i \in[p]\dvtx\vert\mathbf
{y}_i\vert\geq
1/(40\sqrt
{p})\}$.
Since $\|\mathbf{y}\|_1 \leq\sqrt{k}$, the size of $I$ is at most $|
I |
\leq
40\sqrt{kp}$. Now define $\mathbf{y}'\in{\mathbb{R}}^p$ as follows:
$\mathbf
{y}'_i=\mathbf{y}_i$ if
$i\in I$, and $\mathbf{y}'_i=0$ otherwise. By construction, $\mathbf
{y}' \in\hat F$
and $\|\mathbf{y}' - \mathbf{y}\|^2 \leq p \cdot1/(40\sqrt{p})^2 = 1/40^2$.
\end{pf}
We proceed to the discretization of $\hat F$, which uses the following notation.
For $B\subseteq{\mathbb{R}}^p$ and a subset of the coordinates $I
\subseteq[p]$,
let $B_I\subseteq B$ denote the vectors in $B$ whose support is
contained in $I$.
Recall that an \emph{${\varepsilon}$-net} of $B\subseteq{\mathbb{R}}^p$
is a subset $N \subseteq B$ satisfying $B\subseteq N_{\varepsilon}$
and that
for all
$\mathbf{x}\neq\mathbf{y}\in N$, $\Vert\mathbf{x}- \mathbf
{y}\Vert >
{\varepsilon}$.
Setting ${\mathcal{I}}=\{I \subseteq[p]\dvtx |I| = 40\sqrt{pk}\}$,
clearly $\hat F=\bigcup_{I \in{\mathcal{I}}}\hat F_I$.
Let $N_I$ be an ${\varepsilon}$-net of $\hat F_I$ with ${\varepsilon
}=1/40$, and let
$\tilde N$ be the union of all these nets, that is,
\[
\tilde N = \bigcup_{I \in{\mathcal{I}}}
N_I.
\]
Then $\hat F = \bigcup_{I\in{\mathcal{I}}}\hat F_I
\subseteq\bigcup_{I\in{\mathcal{I}}} (N_I)_{1/40}
\subseteq\tilde N_{1/40}
$.
Now using Lemma~\ref{lem:epsNet} and the triangle inequality,
we get that $F \subseteq\hat F_{1/40} \subseteq\tilde N_{1/20}$.
The key to completing the proof is to show that for all sufficiently
large $n$,
%
\begin{equation}
\label{eq:MaxPOverK} \Pr \bigl[\forall\tilde\mathbf{y}\in\tilde N, \tilde\mathbf
{y}^T{\hat{\Sigma}}\tilde \mathbf{y} \leq2(1+\beta) \bigr]
\ge1-e^{-n/10}.
\end{equation}

Before proving this inequality,
let us rely on it to complete the proof of Theorem~\ref{thm:main3}.
Assume the high-probability event in \eqref{eq:MaxPOverK} indeed occurs,
and similarly for Proposition~\ref{prop:LargestEig}, hence
$\lambda_{\max}({\hat{\Sigma}}) \leq(1+{\varepsilon})
(1+\sqrt
{\frac{p}{n}
}+\sqrt{\beta
} )^2$.
Now because $F \subseteq\tilde N_{1/20}$,
for every $\mathbf{y}\in F$ there exists $\tilde\mathbf{y}\in\tilde N$
such that $\mathbf{a}=\mathbf{y}- \tilde\mathbf{y}$ is of length $\|
\mathbf{a}\| \leq1/20$,
and therefore for $Y=\mathbf{y}\mathbf{y}^T$,
%
\begin{equation}
\label{eq:SDPValueOfy} \langle{\hat{\Sigma}}, Y \rangle 
=
\mathbf{y}^T {\hat{\Sigma}}\mathbf{y} = \tilde\mathbf{y}^T
{\hat{\Sigma}}\tilde\mathbf{y}+ 2\mathbf {a}^T {\hat{\Sigma}}\tilde
\mathbf{y}+ \mathbf{a}^T {\hat{\Sigma}}\mathbf{a}.
\end{equation}
The assumption we made using \eqref{eq:MaxPOverK} implies that
$\tilde\mathbf{y}^T {\hat{\Sigma}}\tilde\mathbf{y}\leq2(1+\beta)$.
To bound the two summands, observe that for all $\mathbf{u},\mathbf
{v}\in{\mathbb{R}}^p$ we
have $\mathbf{u}^T {\hat{\Sigma}}\mathbf{v}\leq\Vert\mathbf
{u}\Vert\Vert\mathbf{v}\Vert\lambda_{\max
}({\hat{\Sigma}})$.
Plugging all these into \eqref{eq:SDPValueOfy}, we get
\[
\langle{\hat{\Sigma}}, Y \rangle \leq2(1+\beta) + (1+{\varepsilon}) \biggl(
\frac{2}{20}+\frac
{1}{20^2}\biggr) \biggl(1+\sqrt {\frac{p}{n}}+
\sqrt{\beta} \biggr)^2.
\]
Recall that $\beta\leq\sqrt{\frac{p}{n}}$
and that for sufficiently large $n,p$ we have $\frac{p}{n}\ge20$.
Hence by straightforward manipulations,
we conclude that as $(n,p,k)\to\infty$, with probability tending to one
$\langle{\hat{\Sigma}}, Y \rangle \le\frac{8}{9}\frac{p}{n}$,
which proves Theorem~\ref{thm:main3}.

It remains to prove \eqref{eq:MaxPOverK},
which we do via a union bound argument, using the two lemmas below.
The first lemma estimates the probability that an arbitrary fixed
$\mathbf{y}\in\tilde N$ violates the inequality $\mathbf{y}^T {\hat
{\Sigma}}\mathbf{y}\le
2(1+\beta)$,
and the second one bounds the size of the ${\varepsilon}$-net $\tilde N$.

\begin{lemma} \label{lem:BounOnQuotient}
Under the conditions of Theorem~\ref{thm:main3}, there exists an
integer $n_0 >0$, such that for every $n \ge n_0$ and every $\mathbf
{y}\in{\mathbb{R}}^p$
of length at most $1$
(in particular, every $y\in\tilde N$),
\[
\operatorname{Pr} \bigl[ \mathbf{y}^T {\hat{\Sigma}}\mathbf{y}\geq
2(1+\beta ) \bigr] \leq e^{-n/9}.
\]
\end{lemma}

\begin{pf}
Fix $\mathbf{y}\in{\mathbb{R}}^p$ with $\Vert y\Vert\leq1$, and expand
\[
\mathbf{y}^T {\hat{\Sigma}}\mathbf{y} = \frac{1}n \sum
_{i=1}^n \mathbf{y}^T
\mathbf{x}_i \mathbf {x}_i^T\mathbf{y} =
\frac{1}n \sum_{i=1}^n \langle{
\mathbf{x}_i,\mathbf{y}}\rangle^2 .
\]
Recall from \eqref{eq:defOfDist} that
${\mathbf{x}}_i = \sqrt{\beta}u_i {\mathbf{z}} + \bolds\xi_i$,
where $\bolds\xi_i$ is a vector of independent standard Gaussian
random variables,
and $u_i$ is also a standard Gaussian.
Therefore,
%
\begin{equation}
\label{eq:ithEntry} \langle{\mathbf{x}_i,\mathbf{y}}\rangle= \langle{
\bolds\xi _i,\mathbf{y}}\rangle+ u_i\sqrt{\beta}
\langle{\mathbf {y},\mathbf{z}}\rangle.
\end{equation}
The first term $\langle{\bolds\xi_i,\mathbf{y}}\rangle$
has distribution $N(0,\Vert\mathbf{y}\Vert^2)$.
Since $u_i$ is independent of $\bolds\xi_i$, 
the distribution of $\langle{\mathbf{x}_i,\mathbf{y}}\rangle$ is
just $N(0,\Vert\mathbf{y} \Vert^2+\beta
\langle\mathbf{y},\mathbf{z}\rangle^2)$.
Furthermore, since $\mathbf{y}$ is fixed and the $\bolds\xi_i$'s and
$u_i$'s are all
independent,
the random variables $\langle{\mathbf{x}_i,\mathbf{y}}\rangle$ for
$i=1,\ldots,n$ are i.i.d.,
and thus
\[
\sum_{i=1}^n \langle{
\mathbf{x}_i,\mathbf{y}}\rangle^2 \sim \bigl(\Vert
\mathbf{y}\Vert^2+\beta \langle\mathbf{y} ,\mathbf{z}
\rangle^2 \bigr) \chi^2_n.
\]
Lemma~\ref{lem:XiSquare} with $x=n/9$ implies that $\operatorname
{Pr}[\chi^2_n
\geq
2n] \le e^{-n/9}$. We conclude that with probability at least $1-e^{-n/9}$,
\[
\bigl\langle{\hat{\Sigma}}, \mathbf{y}\mathbf{y}^T \bigr\rangle \leq
\frac{1}{n}\cdot2n\bigl(\Vert\mathbf{y}\Vert^2+\beta\langle
\mathbf {y},\mathbf{z}\rangle^2\bigr) 
\leq2(1+\beta),
\]
where the second inequality uses the Cauchy--Schwarz inequality.
\end{pf}

\begin{lemma} \label{lem:BounOnK}
The ${\varepsilon}$-net $\tilde N$ 
has size $|\tilde N| \leq p^{20\sqrt{pk}}$.
\end{lemma}

\begin{pf}
By the definition of $\tilde N$ and the fact that $|N_I|$ is the same
for all $I$, we can fix arbitrary $I\in{\mathcal{I}}$ and write
%
\begin{equation}
\label{eq:BoundOnK} |\tilde N| \leq\pmatrix{p
\cr
40\sqrt{pk}}|N_I|.
\end{equation}
We thus need to bound $|N_I|$. By definition, $N_I$ is contained in
an axis-aligned subspace of ${\mathbb{R}}^p$ of dimension $p' =
40\sqrt{pk}$,
and we can use the following standard volume argument.
Ignoring henceforth all coordinates outside $I$,
let ${\mathcal{B}}_{r}(\mathbf{x})$ be a closed ball (in ${\mathbb{R}}^{p'}$)
of radius $r>0$
centered at $\mathbf{x}$.
Since $N_I$ is an ${\varepsilon}$-net (of $\hat F_I$),
for every two distinct points in it, $\mathbf{x}\neq\mathbf{y}\in N_I$,
the corresponding balls ${\mathcal{B}}_{{\varepsilon}/2}(\mathbf
{x})$ and
${\mathcal{B}}_{{\varepsilon}/2}(\mathbf{y})$
are disjoint
(as otherwise $\| \mathbf{x}- \mathbf{y}\| \leq{\varepsilon}$).
In addition, the union of these balls ${\mathcal{B}}_{{\varepsilon}
/2}(\mathbf{x})$ over all $\mathbf{x}
\in N_I$
is contained in ${\mathcal{B}}_{1+{\varepsilon}/2}(\mathbf0)$
(because all $\mathbf{x}\in N_I \subseteq\hat F_I$ have length at
most 1).
Recalling that the Euclidean volume of a ball of radius $r>0$ in
dimension $d$
grows with $r$ proportionally to $r^d$,
and plugging in ${\varepsilon}= 1/40$, we obtain
\[
|N_I| \leq\frac{\operatorname{vol}({\mathcal{B}}_{1+{\varepsilon
}/2}(\mathbf
0))}{\operatorname
{vol}({\mathcal{B}}_{{\varepsilon}/2 }(\mathbf0))} \leq \biggl(\frac{1+{\varepsilon}/2}{{\varepsilon}/2}
\biggr)^{p'} = 81^{40\sqrt{pk}}.
\]
Plugging into \eqref{eq:BoundOnK}, we get
$|\tilde N|
\leq (\frac{ep}{40\sqrt{pk}} )^{40\sqrt{pk}} \cdot
81^{40\sqrt{pk}}
\leq p^{20\sqrt{pk}}$.
\end{pf}

Finally, observe that \eqref{eq:MaxPOverK} indeed follows from
Lemmas \ref{lem:BounOnQuotient} and \ref{lem:BounOnK} by a union bound,
\[
\operatorname{Pr} \bigl[\exists\tilde\mathbf{y}\in\tilde N, \tilde
\mathbf{y}^T {\hat{\Sigma}}\tilde \mathbf{y} \geq2(1+\beta) \bigr]
\leq p^{20\sqrt{pk}}\cdot e^{-n/9} \leq e^{-n/10},
\]
where the last inequality follows from the assumption
in Theorem~\ref{thm:main3} that $k/( p/\log^2p) \to0$ and that
$p/n\to c$.
This completes the proof of \eqref{eq:MaxPOverK}
and of Theorem~\ref{thm:main3}.

\section{Deferred proofs from Section 
\texorpdfstring{\protect\ref{sec:prelims}}{5}
(preliminaries)}
\label{sec:PrelimProofs}

\mbox{}
\begin{pf*}{Proof of Proposition~\ref{prop:LargestEig}}
Let us rotate ${\mathbb{R}}^p$ so that the spike $\mathbf{z}$ becomes
the first standard
basis vector $\mathbf{e}_1$.
Obviously, $\lambda_{\max}({\hat{\Sigma}})$ would not change at all,
and since the normal distribution is rotation invariant,
the noise would still be normally distributed.
In effect, we may assume henceforth that $\mathbf{z}=\mathbf{e}_1$.
Recalling from \eqref{eq:defOfDist} that the samples are given by
$\mathbf{x}
_i=\sqrt{\beta}u_i {\mathbf{z}} + \bolds\xi_i$,
we can write ${\hat{\Sigma}}= \frac{1}{n}\sum_{i=1}^n \mathbf{x}_i
\mathbf{x}_i^T$ as
\[
{\hat{\Sigma}}= \frac{1}{n} (\sqrt\beta U + \Xi ) (\sqrt \beta
U + \Xi )^T,
\]
where $U$ is a $p \times n$ matrix whose first row is $(u_1,\ldots
,u_n)$ and the remaining rows are zero (recall $\mathbf{z}=\mathbf{e}_1$),
and $\Xi$ is an $p \times n$ matrix whose $i$th column is $\bolds\xi_i$.
Let $\Vert A\Vert=\sqrt{\lambda_{\max}(A^TA)}$ be the spectral norm
of a
matrix $A$.
Using also $\Vert A\Vert=\Vert A^T\Vert$ and the triangle inequality,
%
\begin{equation}
\label{eq:EigEq1} \lambda_{\max}({\hat{\Sigma}}) 
=
\frac{1}{n} \|\sqrt\beta U + \Xi\|^2 \le\frac{1}{n} \bigl(
\sqrt{\beta} \Vert U\Vert + \Vert\Xi\Vert \bigr)^2.
\end{equation}

The matrix $\Xi^T\Xi$ follows a Wishart distribution (note that the
roles of $p$ and $n$ are reversed). Therefore by \cite{ElKaroui08},
Theorem~2, which applies to the regime $p/n \to0$ and $p/n \to
\infty$, and by \cite{Johnstone01}, Theorem~1.1, which applies to $p/n
\to c \in(0,\infty)$,
we know that with probability tending to one, 
\[
\|\Xi\|^2 = \lambda_{\max}\bigl(\Xi^T
\Xi\bigr) \le(1+{\varepsilon}_1) (\sqrt {p}+\sqrt {n}
)^2,
\]
for some ${\varepsilon}_1={\varepsilon}_1(n) \to0$. Since $U^TU$ has
rank one,
$\|U\|^2=\lambda_{\max}(U^TU)=\operatorname{tr}(U^TU)=\sum_{i=1}^n
u_i^2 \sim
\chi
^2_n $. Lemma~\ref{lem:XiSquare} with $x=3 \log n$ implies that with
probability at least $1-1/n^3$,
\[
\Vert U\Vert^2 \le \biggl(1+O \biggl(\frac{\log n} {n}
\biggr) \biggr)n,
\]
and thus with probability tending to one,
$\Vert U\Vert \le(1+{\varepsilon}_2)\sqrt{n}$ for some
${\varepsilon}_2={\varepsilon}_2(n) \to0$.

Plugging these bounds 
into \eqref{eq:EigEq1},
we conclude that with probability tending to one as $n \to\infty$,
\begin{eqnarray*}
\lambda_{\max}({\hat{\Sigma}}) 
&\le&
\biggl[ (1+{\varepsilon}_2)\sqrt{\beta} + (1+{\varepsilon
}_1) \biggl(1+\sqrt {\frac{p}{n} } \biggr)
\biggr]^2
\\
&\le& \biggl[ (1+{\varepsilon}_1+{\varepsilon}_2)
\biggl(\sqrt{\beta } + 1+\sqrt {\frac{p}{n} } \biggr) \biggr]^2
,
\end{eqnarray*}
which completes the proof of Proposition~\ref{prop:LargestEig}.
\end{pf*}


\begin{pf*}{Proof of Proposition~\ref{fact:ledoit}}
Starting with $\operatorname{tr}({\hat{\Sigma}})$,
observe that ${\hat{\Sigma}}_{ii}\sim\frac{1}{n}\chi^2_n$.
Lemma~\ref{lem:XiSquare} with $x=5\ln p$ implies that with probability
at least $1-1/p^5$, $\chi^2_n \leq(1+{\varepsilon}_1)n$ for
${\varepsilon}_1 =
O(\sqrt
{(\log p)/n}) \to0$.
Taking a union bound over $i=1,\ldots,p$, we obtain that with
probability at least $1-1/p^4$, all entries ${\hat{\Sigma}}_{ii} \in
[1-{\varepsilon}
_1,1+{\varepsilon}_1]$,
which implies
\[
\operatorname{tr}({\hat{\Sigma}})=\sum_{i=1}^p
{\hat{\Sigma }}_{ii} \in\bigl[(1-{\varepsilon}_1)p,(1+{
\varepsilon}_1)p\bigr].
\]

We now turn to bound
\[
\operatorname{tr}\bigl({\hat{\Sigma}}^2\bigr) = \sum
_{i,j=1}^p {\hat{\Sigma}}_{ij}^2
= \sum_{i=1}^p {\hat{\Sigma}}_{ii}^2
+ \sum_{i=1}^p \sum
_{j \ne i} {\hat{\Sigma}}_{ij}^2.
\]
By the preceding paragraph, with probability at least $1-1/p^4$,
$\sum_{i=1}^p {\hat{\Sigma}}_{ii}^2 \in[(1-{\varepsilon
}_1)^2p,(1+{\varepsilon}_1)^2p]$.
Using the notation of \eqref{eq:defOfDist},
we write off-diagonal entries in ${\hat{\Sigma}}$ as
${\hat{\Sigma}}_{ij}= \frac{1}{n}\sum_{s=1}^n \bolds\xi
_{si}\bolds\xi_{sj}:=\frac{1}{n}
\bolds\rho_i^T\bolds\rho_j$,
where $\bolds\rho_i = (\bolds\xi_{si})_{s=1}^n$,
and notice that $\bolds\rho_1,\ldots,\bolds\rho_p$ are independent.

Now fix $i$ and condition on $\bolds\rho_i$.
Then Lemma~\ref{cl:CLT} implies that each off-diagonal entry along row $i$
is distributed ${\hat{\Sigma}}_{ij} \sim\frac{1}{n} \Vert\bolds
\rho _i\Vert\cdot
\hat
y_j$, $\hat y_j \sim N(0,1$).
Moreover the $\hat y_j$'s (for different $j\neq i$) are independent, hence,
$\sum_{j \ne i} {\hat{\Sigma}}_{ij}^2 \sim\frac{1}{n^2} \Vert
\bolds\rho_i\Vert^2
\chi
_{p-1}^2$.
Using Lemma~\ref{lem:XiSquare} with $x = 4 \log p$,
with probability at least $1-1/p^4$,
\[
\sum_{j \ne i} {\hat{\Sigma}}_{ij}^2
\in \biggl[(1 -{\varepsilon}_2) (p-1)\cdot\frac{1}{n^2} \Vert
\bolds \rho_i\Vert^2 ,(1 +{\varepsilon}_2)
(p-1)\cdot\frac{1}{n^2} \Vert\bolds\rho _i\Vert^2
\biggr],
\]
for ${\varepsilon}_2 = O(\sqrt{(\log p)/ p})$.

Next, remove the conditioning on $\bolds\rho_i$ (still for a fixed $i$),
observing that $\|\bolds\rho_i\|^2 \sim\chi^2_n$.
Lemma~\ref{lem:XiSquare} with $x = 4 \log p$ then implies that
with probability at least $1-1/p^4$,
we have $\|\bolds\rho_i\|^2 \in[(1 -{\varepsilon}_3)n,(1
+{\varepsilon}_3)n]$
for ${\varepsilon}_3 = O(\sqrt{(\log p)/n})$.

Finally, taking the union bound over rows $i=1,\ldots,p$
and also the sum along the diagonal,
with probability at least $1-3/p^3$,
\[
\operatorname{tr}\bigl({\hat{\Sigma}}^2\bigr) \le(1+{
\varepsilon}_1)^2p + (1+{\varepsilon}_2)
(p-1)\cdot\frac
{1}{n^2}\cdot(1+{\varepsilon}_3)n \le(1+{
\varepsilon}_4)p \biggl(1+\frac{p}{n} \biggr),
\]
for a suitably chosen ${\varepsilon}_4={\varepsilon}_4(n)\to0$.
Similarly, $\operatorname{tr}({\hat{\Sigma}}^2) \ge(1-{\varepsilon}
_5)p (1+\frac{p}{n} )$
for ${\varepsilon}_5={\varepsilon}_5(n)\to0$.
To complete the proof of Proposition~\ref{fact:ledoit},
set ${\varepsilon}=\max\{{\varepsilon}_1,{\varepsilon
}_4,{\varepsilon}_5\}$.
\end{pf*}

%
%

%



\printaddresses

\begin{thebibliography}{42}

\bibitem{AKSClique}
%
\begin{barticle}[mr]
\bauthor{\bsnm{Alon},~\bfnm{Noga}\binits{N.}},
\bauthor{\bsnm{Krivelevich},~\bfnm{Michael}\binits{M.}} \AND
\bauthor{\bsnm{Sudakov},~\bfnm{Benny}\binits{B.}}
(\byear{1998}).
\btitle{Finding a large hidden clique in a~random graph}.
\bjournal{{R}andom {S}tructures {A}lgorithms}
\bvolume{13}
\bpages{457--466}.
\bid{doi={10.1002/(SICI)1098-2418(199810/12)13:3/4<457::AID-RSA14>3.3.CO;2-K},
issn={1042-9832}, mr={1662795}}\vadjust{\goodbreak}
\end{barticle}
%

\bptok{imsref}%
\endbibitem

\bibitem{Ames:2011}
%
\begin{barticle}[mr]
\bauthor{\bsnm{Ames},~\bfnm{Brendan~P.~W.}\binits{B.~P.~W.}} \AND
\bauthor{\bsnm{Vavasis},~\bfnm{Stephen~A.}\binits{S.~A.}}
(\byear{2011}).
\btitle{Nuclear norm minimization for the planted clique and biclique
problems}.
\bjournal{Math. Program.}
\bvolume{129}
\bpages{69--89}.
\bid{doi={10.1007/s10107-011-0459-x}, issn={0025-5610}, mr={2831403}}
\end{barticle}
%

\bptok{imsref}%
\endbibitem

\bibitem{AminiWain09}
%
\begin{barticle}[mr]
\bauthor{\bsnm{Amini},~\bfnm{Arash~A.}\binits{A.~A.}} \AND
\bauthor{\bsnm{Wainwright},~\bfnm{Martin~J.}\binits{M.~J.}}
(\byear{2009}).
\btitle{High-dimensional analysis of semidefinite relaxations for
sparse principal components}.
\bjournal{Ann. Statist.}
\bvolume{37}
\bpages{2877--2921}.
\bid{doi={10.1214/08-AOS664}, issn={0090-5364}, mr={2541450}}
\end{barticle}
%

\bptok{imsref}%
\endbibitem

\bibitem{Anderson84}
%
\begin{bbook}[mr]
\bauthor{\bsnm{Anderson},~\bfnm{T.~W.}\binits{T.~W.}}
(\byear{1984}).
\btitle{An Introduction to Multivariate Statistical Analysis},
\bedition{2nd} ed.
\bpublisher{Wiley},
\blocation{New York}.
\bid{mr={0771294}}
\end{bbook}
%

\bptok{imsref}%
\endbibitem

\bibitem{BaikSilverstein06}
%
\begin{barticle}[mr]
\bauthor{\bsnm{Baik},~\bfnm{Jinho}\binits{J.}} \AND
\bauthor{\bsnm{Silverstein},~\bfnm{Jack~W.}\binits{J.~W.}}
(\byear{2006}).
\btitle{Eigenvalues of large sample covariance matrices of spiked
population models}.
\bjournal{J. Multivariate Anal.}
\bvolume{97}
\bpages{1382--1408}.
\bid{doi={10.1016/j.jmva.2005.08.003}, issn={0047-259X}, mr={2279680}}
\end{barticle}
%

\bptok{imsref}%
\endbibitem

\bibitem{RigolletClique}
%
\begin{barticle}[mr]
\bauthor{\bsnm{Berthet},~\bfnm{Quentin}\binits{Q.}} \AND
\bauthor{\bsnm{Rigollet},~\bfnm{Philippe}\binits{P.}}
(\byear{2013}).
\btitle{Optimal detection of sparse principal components in high dimension}.
\bjournal{Ann. Statist.}
\bvolume{41}
\bpages{1780--1815}.
\bid{doi={10.1214/13-AOS1127}, issn={0090-5364}, mr={3127849}}
\end{barticle}
%

\bptok{imsref}%
\endbibitem

\bibitem{RigolletCliqueCOLT}
%
\begin{bmisc}[auto]
\bauthor{\bsnm{Berthet},~\bfnm{Quentin}\binits{Q.}} \AND
\bauthor{\bsnm{Rigollet},~\bfnm{Philippe}\binits{P.}}
(\byear{2013}).
\bhowpublished{Complexity theoretic lower bounds for sparse principal component
detection.
In \textit{COLT} 1046--1066. \surl{JMLR.org}.}
\end{bmisc}
%

\bptok{imsref}%
\endbibitem

\bibitem{BickelLevina06}
%
\begin{barticle}[mr]
\bauthor{\bsnm{Bickel},~\bfnm{Peter~J.}\binits{P.~J.}} \AND
\bauthor{\bsnm{Levina},~\bfnm{Elizaveta}\binits{E.}}
(\byear{2008}).
\btitle{Regularized estimation of large covariance matrices}.
\bjournal{Ann. Statist.}
\bvolume{36}
\bpages{199--227}.
\bid{doi={10.1214/009053607000000758}, issn={0090-5364}, mr={2387969}}
\end{barticle}
%

\bptok{imsref}%
\endbibitem

\bibitem{BJNP}
%
\begin{barticle}[mr]
\bauthor{\bsnm{Birnbaum},~\bfnm{Aharon}\binits{A.}},
\bauthor{\bsnm{Johnstone},~\bfnm{Iain~M.}\binits{I.~M.}},
\bauthor{\bsnm{Nadler},~\bfnm{Boaz}\binits{B.}} \AND
\bauthor{\bsnm{Paul},~\bfnm{Debashis}\binits{D.}}
(\byear{2013}).
\btitle{Minimax bounds for sparse PCA with noisy high-dimensional data}.
\bjournal{Ann. Statist.}
\bvolume{41}
\bpages{1055--1084}.
\bid{doi={10.1214/12-AOS1014}, issn={0090-5364}, mr={3113803}}
\end{barticle}
%

\bptok{imsref}%
\endbibitem

\bibitem{CaiMaWu13}
%
\begin{barticle}[mr]
\bauthor{\bsnm{Cai},~\bfnm{T.~Tony}\binits{T.~T.}},
\bauthor{\bsnm{Ma},~\bfnm{Zongming}\binits{Z.}} \AND
\bauthor{\bsnm{Wu},~\bfnm{Yihong}\binits{Y.}}
(\byear{2013}).
\btitle{Sparse PCA: Optimal rates and adaptive estimation}.
\bjournal{Ann. Statist.}
\bvolume{41}
\bpages{3074--3110}.
\bid{doi={10.1214/13-AOS1178}, issn={0090-5364}, mr={3161458}}
\end{barticle}
%

\bptok{imsref}%
\endbibitem

\bibitem{AspremontBannerjeeGhaoui07}
%
\begin{barticle}[mr]
\bauthor{\bsnm{d'Aspremont},~\bfnm{Alexandre}\binits{A.}},
\bauthor{\bsnm{Banerjee},~\bfnm{Onureena}\binits{O.}} \AND
\bauthor{\bsnm{El Ghaoui},~\bfnm{Laurent}\binits{L.}}
(\byear{2008}).
\btitle{First-order methods for sparse covariance selection}.
\bjournal{SIAM J. Matrix Anal. Appl.}
\bvolume{30}
\bpages{56--66}.
\bid{doi={10.1137/060670985}, issn={0895-4798}, mr={2399568}}
\end{barticle}
%

\bptok{imsref}%
\endbibitem

\bibitem{AspremontEtAlSDP07}
%
\begin{barticle}[auto:parserefs-M02]
\bauthor{\bsnm{d'Aspremont},~\bfnm{A.}\binits{A.}},
\bauthor{\bsnm{El-Ghaoui},~\bfnm{L.}\binits{L.}},
\bauthor{\bsnm{Jordan},~\bfnm{M.}\binits{M.}} \AND
\bauthor{\bsnm{Lanckriet},~\bfnm{G.}\binits{G.}}
(\byear{2004}).
\btitle{A direct formulation for sparse PCA using semidefinite programming}.
\bjournal{SIAM Rev.}
\bvolume{49}
\bpages{434--448}.
\bid{doi={10.1137/050645506}}
\end{barticle}
%

\bptok{imsref}%
\endbibitem

\bibitem{DPaul07}
%
\begin{barticle}[mr]
\bauthor{\bsnm{Debashis},~\bfnm{Paul}\binits{P.}}
(\byear{2007}).
\btitle{Asymptotics of sample eigenstructure for a large dimensional
spiked covariance model}.
\bjournal{Statist. Sinica}
\bvolume{17}
\bpages{1617--1642}.
\bid{issn={1017-0405}, mr={2399865}}
\end{barticle}
%

\bptok{imsref}%
\endbibitem

\bibitem{PeresClique}
%
\begin{barticle}[mr]
\bauthor{\bsnm{Dekel},~\bfnm{Yael}\binits{Y.}},
\bauthor{\bsnm{Gurel-Gurevich},~\bfnm{Ori}\binits{O.}} \AND
\bauthor{\bsnm{Peres},~\bfnm{Yuval}\binits{Y.}}
(\byear{2014}).
\btitle{Finding hidden cliques in linear time with high probability}.
\bjournal{Combin. Probab. Comput.}
\bvolume{23}
\bpages{29--49}.
\bid{doi={10.1017/S096354831300045X}, issn={0963-5483}, mr={3197965}}
\end{barticle}
%

\bptok{imsref}%
\endbibitem

\bibitem{MontanariClique}
%
\begin{bmisc}[auto:parserefs-M02]
\bauthor{\bsnm{Deshpande},~\bfnm{Y.}\binits{Y.}} \AND
\bauthor{\bsnm{Montanari},~\bfnm{A.}\binits{A.}}
(\byear{2013}).
\bhowpublished{Finding hidden cliques of size {$\sqrt{n/e}$} in nearly
linear time.
Available at \arxivurl{arXiv:1304.7047}.}
\end{bmisc}
%

\bptok{imsref}%
\endbibitem

\bibitem{MontanariCovar}
%
\begin{bmisc}[auto:parserefs-M02]
\bauthor{\bsnm{Deshpande},~\bfnm{Y.}\binits{Y.}} \AND
\bauthor{\bsnm{Montanari},~\bfnm{A.}\binits{A.}}
(\byear{2014}).
\bhowpublished{Information-theoretically optimal sparse PCA.
Available at \arxivurl{arXiv:1402.2238}}.
\end{bmisc}
%

\bptok{imsref}%
\endbibitem

\bibitem{ElKaroui08}
%
\begin{bmisc}[auto:parserefs-M02]
\bauthor{\bsnm{El-Karoui},~\bfnm{N.}\binits{N.}}
(\byear{2003}).
\bhowpublished{On the largest eigenvalue of wishart matrices with identity
covariance when {$n,p$} and {$p/n \to\infty$}.
Available at \arxivurl{arXiv:math/0309355}}.
\end{bmisc}
%

\bptok{imsref}%
\endbibitem

\bibitem{FK00}
%
\begin{barticle}[mr]
\bauthor{\bsnm{Feige},~\bfnm{Uriel}\binits{U.}} \AND
\bauthor{\bsnm{Krauthgamer},~\bfnm{Robert}\binits{R.}}
(\byear{2000}).
\btitle{Finding and certifying a large hidden clique in a semirandom graph}.
\bjournal{Random Structures Algorithms}
\bvolume{16}
\bpages{195--208}.
\bid
{doi={10.1002/(SICI)1098-2418(200003)16:2<195::AID-RSA5>3.3.CO;2-1},
issn={1042-9832}, mr={1742351}}
\end{barticle}
%

\bptok{imsref}%
\endbibitem

\bibitem{FeigeClique}
%
\begin{bincollection}[mr]
\bauthor{\bsnm{Feige},~\bfnm{Uriel}\binits{U.}} \AND
\bauthor{\bsnm{Ron},~\bfnm{Dorit}\binits{D.}}
(\byear{2010}).
\btitle{Finding hidden cliques in linear time}.
In \bbooktitle{21st {I}nternational {M}eeting on {P}robabilistic,
{C}ombinatorial, and {A}symptotic {M}ethods in the {A}nalysis of
{A}lgorithms ({A}of{A}'10)}
\bpages{189--203}.
\bpublisher{Assoc. Discrete Math. Theor. Comput. Sci.},
\blocation{Nancy}.
\bid{mr={2735341}}
\end{bincollection}
%

\bptok{imsref}%
\endbibitem

\bibitem{Johnstone01}
%
\begin{barticle}[mr]
\bauthor{\bsnm{Johnstone},~\bfnm{Iain~M.}\binits{I.~M.}}
(\byear{2001}).
\btitle{On the distribution of the largest eigenvalue in principal
components analysis}.
\bjournal{Ann. Statist.}
\bvolume{29}
\bpages{295--327}.
\bid{doi={10.1214/aos/1009210544}, issn={0090-5364}, mr={1863961}}
\end{barticle}
%

\bptok{imsref}%
\endbibitem

\bibitem{Johnstone.Lu2009Consistency}
%
\begin{barticle}[mr]
\bauthor{\bsnm{Johnstone},~\bfnm{Iain~M.}\binits{I.~M.}} \AND
\bauthor{\bsnm{Lu},~\bfnm{Arthur~Yu}\binits{A.~Y.}}
(\byear{2009}).
\btitle{On consistency and sparsity for principal components analysis
in high dimensions}.
\bjournal{J. Amer. Statist. Assoc.}
\bvolume{104}
\bpages{682--693}.
\bid{doi={10.1198/jasa.2009.0121}, issn={0162-1459}, mr={2751448}}
\end{barticle}
%

\bptok{imsref}%
\endbibitem

\bibitem{JoliffePCA}
%
\begin{bbook}[mr]
\bauthor{\bsnm{Jolliffe},~\bfnm{I.~T.}\binits{I.~T.}}
(\byear{2002}).
\btitle{Principal Component Analysis},
\bedition{2nd} ed.
\bpublisher{Springer},
\blocation{New York}.
\bid{mr={2036084}}
\end{bbook}
%

\bptok{imsref}%
\endbibitem

\bibitem{JolliffeTrendafilovUddin03}
%
\begin{barticle}[mr]
\bauthor{\bsnm{Jolliffe},~\bfnm{Ian~T.}\binits{I.~T.}},
\bauthor{\bsnm{Trendafilov},~\bfnm{Nickolay~T.}\binits{N.~T.}} \AND
\bauthor{\bsnm{Uddin},~\bfnm{Mudassir}\binits{M.}}
(\byear{2003}).
\btitle{A modified principal component technique based on the {LASSO}}.
\bjournal{J. Comput. Graph. Statist.}
\bvolume{12}
\bpages{531--547}.
\bid{doi={10.1198/1061860032148}, issn={1061-8600}, mr={2002634}}
\end{barticle}
%

\bptok{imsref}%
\endbibitem

\bibitem{LAURENT-MASSART90}
%
\begin{barticle}[mr]
\bauthor{\bsnm{Laurent},~\bfnm{B.}\binits{B.}} \AND
\bauthor{\bsnm{Massart},~\bfnm{P.}\binits{P.}}
(\byear{2000}).
\btitle{Adaptive estimation of a quadratic functional by model selection}.
\bjournal{Ann. Statist.}
\bvolume{28}
\bpages{1302--1338}.
\bid{doi={10.1214/aos/1015957395}, issn={0090-5364}, mr={1805785}}
\end{barticle}
%

\bptok{imsref}%
\endbibitem

\bibitem{Ledoit02}
%
\begin{barticle}[mr]
\bauthor{\bsnm{Ledoit},~\bfnm{Olivier}\binits{O.}} \AND
\bauthor{\bsnm{Wolf},~\bfnm{Michael}\binits{M.}}
(\byear{2002}).
\btitle{Some hypothesis tests for the covariance matrix when the
dimension is large compared to the sample size}.
\bjournal{Ann. Statist.}
\bvolume{30}
\bpages{1081--1102}.
\bid{doi={10.1214/aos/1031689018}, issn={0090-5364}, mr={1926169}}
\end{barticle}
%

\bptok{imsref}%
\endbibitem

\bibitem{LeiVu14}
%
\begin{barticle}[mr]
\bauthor{\bsnm{Lei},~\bfnm{Jing}\binits{J.}} \AND
\bauthor{\bsnm{Vu},~\bfnm{Vincent~Q.}\binits{V.~Q.}}
(\byear{2015}).
\btitle{Sparsistency and agnostic inference in sparse {PCA}}.
\bjournal{Ann. Statist.}
\bvolume{43}
\bpages{299--322}.
\bid{doi={10.1214/14-AOS1273}, issn={0090-5364}, mr={3311861}}
\bptnote{check year}%
\end{barticle}
%

\bptok{imsref}%
\endbibitem

\bibitem{ZhaosongZhang:2012}
%
\begin{barticle}[mr]
\bauthor{\bsnm{Lu},~\bfnm{Zhaosong}\binits{Z.}} \AND
\bauthor{\bsnm{Zhang},~\bfnm{Yong}\binits{Y.}}
(\byear{2012}).
\btitle{An augmented {L}agrangian approach for sparse principal
component analysis}.
\bjournal{Math. Program.}
\bvolume{135}
\bpages{149--193}.
\bid{doi={10.1007/s10107-011-0452-4}, issn={0025-5610}, mr={2968253}}
\end{barticle}
%

\bptok{imsref}%
\endbibitem

\bibitem{Ma2013}
%
\begin{barticle}[mr]
\bauthor{\bsnm{Ma},~\bfnm{Zongming}\binits{Z.}}
(\byear{2013}).
\btitle{Sparse principal component analysis and iterative thresholding}.
\bjournal{Ann. Statist.}
\bvolume{41}
\bpages{772--801}.
\bid{doi={10.1214/13-AOS1097}, issn={0090-5364}, mr={3099121}}
\end{barticle}
%

\bptok{imsref}%
\endbibitem

\bibitem{Moghaddam06generalizedspectral}
%
\begin{binproceedings}[auto:parserefs-M02]
\bauthor{\bsnm{Moghaddam},~\bfnm{B.}\binits{B.}},
\bauthor{\bsnm{Weiss},~\bfnm{S.}\binits{S.}} \AND
\bauthor{\bsnm{Avidan},~\bfnm{Y.}\binits{Y.}}
(\byear{2006}).
\btitle{Generalized spectral bounds for sparse {LDA}}.
In \bbooktitle{Proceedings of the 23rd International Conference on
Machine Learning}
\bpages{641--648}.
\bpublisher{ACM},
\blocation{New York}.
\bid{doi={10.1145/1143844.1143925}}
\end{binproceedings}
%

\bptok{imsref}%
\endbibitem

\bibitem{Moghaddam06spectralbounds}
%
\begin{bincollection}[auto:parserefs-M02]
\bauthor{\bsnm{Moghaddam},~\bfnm{B.}\binits{B.}},
\bauthor{\bsnm{Weiss},~\bfnm{Y.}\binits{Y.}} \AND
\bauthor{\bsnm{Avidan},~\bfnm{S.}\binits{S.}}
(\byear{2006}).
\btitle{Spectral bounds for sparse PCA: Exact and greedy algorithms}.
In \bbooktitle{Advances in Neural Information Processing Systems}
\bpages{915--922}.
\bpublisher{MIT Press},
\blocation{Cambridge}.
\end{bincollection}
%

\bptok{imsref}%
\endbibitem

\bibitem{Muirhead1982}
%
\begin{bbook}[mr]
\bauthor{\bsnm{Muirhead},~\bfnm{Robb~J.}\binits{R.~J.}}
(\byear{1982}).
\btitle{Aspects of Multivariate Statistical Theory}.
\bpublisher{Wiley},
\blocation{New York}.
\bid{mr={0652932}}
\end{bbook}
%

\bptok{imsref}%
\endbibitem

\bibitem{Nadler08}
%
\begin{barticle}[mr]
\bauthor{\bsnm{Nadler},~\bfnm{Boaz}\binits{B.}}
(\byear{2008}).
\btitle{Finite sample approximation results for principal component
analysis: A~matrix perturbation approach}.
\bjournal{Ann. Statist.}
\bvolume{36}
\bpages{2791--2817}.
\bid{doi={10.1214/08-AOS618}, issn={0090-5364}, mr={2485013}}
\end{barticle}
%

\bptok{imsref}%
\endbibitem

\bibitem{Natarajan1995}
%
\begin{barticle}[mr]
\bauthor{\bsnm{Natarajan},~\bfnm{B.~K.}\binits{B.~K.}}
(\byear{1995}).
\btitle{Sparse approximate solutions to linear systems}.
\bjournal{SIAM J. Comput.}
\bvolume{24}
\bpages{227--234}.
\bid{doi={10.1137/S0097539792240406}, issn={0097-5397}, mr={1320206}}
\end{barticle}
%

\bptok{imsref}%
\endbibitem

\bibitem{Shen:2008}
%
\begin{barticle}[mr]
\bauthor{\bsnm{Shen},~\bfnm{Haipeng}\binits{H.}} \AND
\bauthor{\bsnm{Huang},~\bfnm{Jianhua~Z.}\binits{J.~Z.}}
(\byear{2008}).
\btitle{Sparse principal component analysis via regularized low rank
matrix approximation}.
\bjournal{J. Multivariate Anal.}
\bvolume{99}
\bpages{1015--1034}.
\bid{doi={10.1016/j.jmva.2007.06.007}, issn={0047-259X}, mr={2419336}}
\end{barticle}
%

\bptok{imsref}%
\endbibitem

\bibitem{STEWART90}
%
\begin{bbook}[mr]
\bauthor{\bsnm{Stewart},~\bfnm{G.~W.}\binits{G.~W.}} \AND
\bauthor{\bsnm{Sun},~\bfnm{Ji~Guang}\binits{J.~G.}}
(\byear{1990}).
\btitle{Matrix Perturbation Theory}.
\bpublisher{Academic Press},
\blocation{Boston, MA}.
\bid{mr={1061154}}
\end{bbook}
%

\bptok{imsref}%
\endbibitem

\bibitem{sedumi}
%
\begin{barticle}[mr]
\bauthor{\bsnm{Sturm},~\bfnm{Jos~F.}\binits{J.~F.}}
(\byear{1999}).
\btitle{Using {S}e{D}u{M}i 1.02, a MATLAB toolbox for optimization
over symmetric cones}.
\bjournal{Optim. Methods Softw.}
\bvolume{11/12}
\bpages{625--653}.
\bid{doi={10.1080/10556789908805766}, issn={1055-6788}, mr={1778433}}
\end{barticle}
%

\bptok{imsref}%
\endbibitem

\bibitem{LeiVu2013}
%
\begin{barticle}[mr]
\bauthor{\bsnm{Vu},~\bfnm{Vincent~Q.}\binits{V.~Q.}} \AND
\bauthor{\bsnm{Lei},~\bfnm{Jing}\binits{J.}}
(\byear{2013}).
\btitle{Minimax sparse principal subspace estimation in high dimensions}.
\bjournal{Ann. Statist.}
\bvolume{41}
\bpages{2905--2947}.
\bid{doi={10.1214/13-AOS1151}, issn={0090-5364}, mr={3161452}}
\end{barticle}
%

\bptok{imsref}%
\endbibitem

\bibitem{Samworth14}
%
\begin{bmisc}[auto:parserefs-M02]
\bauthor{\bsnm{Wang},~\bfnm{T.}\binits{T.}},
\bauthor{\bsnm{Berthet},~\bfnm{Q.}\binits{Q.}} \AND
\bauthor{\bsnm{Samworth},~\bfnm{R.}\binits{R.}}
(\byear{2014}).
\bhowpublished{Statistical and computational trade-offs in estimation of
sparse principal components.
Available at \arxivurl{arXiv:1408.5369}.}
\end{bmisc}
%

\bptok{imsref}%
\endbibitem

\bibitem{WittenTibshirani:2009}
%
\begin{barticle}[auto:parserefs-M02]
\bauthor{\bsnm{Witten},~\bfnm{D.}\binits{D.}},
\bauthor{\bsnm{Tibshirani},~\bfnm{R.}\binits{R.}} \AND
\bauthor{\bsnm{Hastie},~\bfnm{R.}\binits{R.}}
(\byear{2009}).
\btitle{A penalized matrix decomposition, with applications to sparse
principal components and canonical correlation analysis}.
\bjournal{Biostatistics}
\bvolume{10}
\bpages{515--534}.
\end{barticle}
%

\bptok{imsref}%
\endbibitem

\bibitem{ZhangZhaSimon}
%
\begin{barticle}[mr]
\bauthor{\bsnm{Zhang},~\bfnm{Zhenyue}\binits{Z.}},
\bauthor{\bsnm{Zha},~\bfnm{Hongyuan}\binits{H.}} \AND
\bauthor{\bsnm{Simon},~\bfnm{Horst}\binits{H.}}
(\byear{2002}).
\btitle{Low-rank approximations with sparse factors. I. {B}asic
algorithms and error analysis}.
\bjournal{SIAM J. Matrix Anal. Appl.}
\bvolume{23}
\bpages{706--727 (electronic)}.
\bid{doi={10.1137/S0895479899359631}, issn={0895-4798}, mr={1896815}}
\bptnote{check year}%
\end{barticle}
%

\bptok{imsref}%
\endbibitem

\bibitem{WangLuLiu14}
%
\begin{bmisc}[auto:parserefs-M02]
\bauthor{\bsnm{Zhaoran Wang},~\bfnm{Z.}\binits{Z.}},
\bauthor{\bsnm{Lu},~\bfnm{H.}\binits{H.}} \AND
\bauthor{\bsnm{Liu},~\bfnm{H.}\binits{H.}}
(\byear{2014}).
\bhowpublished{Nonconvex statistical optimization: Minimax-optimal sparse pca
in polynomial time.
Available at \arxivurl{arXiv:1408.5352}.}
\end{bmisc}
%

\bptok{imsref}%
\endbibitem

\bibitem{ZouHastieTibshirani06}
%
\begin{barticle}[mr]
\bauthor{\bsnm{Zou},~\bfnm{Hui}\binits{H.}},
\bauthor{\bsnm{Hastie},~\bfnm{Trevor}\binits{T.}} \AND
\bauthor{\bsnm{Tibshirani},~\bfnm{Robert}\binits{R.}}
(\byear{2006}).
\btitle{Sparse principal component analysis}.
\bjournal{J.~Comput. Graph. Statist.}
\bvolume{15}
\bpages{265--286}.
\bid{doi={10.1198/106186006X113430}, issn={1061-8600}, mr={2252527}}
\end{barticle}
%

\bptok{imsref}%
\endbibitem

\end{thebibliography}
\end{document}